\documentclass[11pt,leqno,twoside]{article}

\usepackage{amsfonts,amsmath,amsthm,amssymb}
\usepackage{enumerate}
\usepackage{graphics,graphicx,subfigure}

\usepackage{color}
\usepackage{todonotes}
\usepackage{cancel}
\usepackage{url}
\usepackage{hyperref}
\usepackage{makeidx}
\usepackage{showidx}
\usepackage{multicol}        % used for the two-column index
\usepackage{xspace}
\usepackage{stmaryrd}        % for \llbracket, \rrbracket
\usepackage{pifont}          % for \ding{227}
\usepackage{fancybox}        % for oval boxes
\usepackage{bm}

 \usepackage{fancyhdr}
\theoremstyle{plain}
 \theoremstyle{definition}
 \newtheorem{lem}{Lemma}
 \newtheorem{defn}[lem]{Definition}
 \newtheorem{thm}[lem]{Theorem}
 \newtheorem{prop}[lem]{Proposition}
 \newtheorem{cor}[lem]{Corollary}
 \newtheorem{notn}[lem]{Notations}
 \newtheorem{pb}[lem]{Problem}
 \newtheorem{form}[lem]{Formulae}
 
 \newtheorem*{rk}{Remark}
 \newtheorem*{com}{Comment}
 \newtheorem*{ex}{Example}
 \theoremstyle{remark}

 \newcommand{\blem}{\begin{lem}}
 \newcommand{\elem}{\end{lem}}
 \newcommand{\bdefn}{\begin{defn}}
 \newcommand{\edefn}{\end{defn}}
 \newcommand{\bthm}{\begin{thm} }
 \newcommand{\ethm}{\end{thm}}
 \newcommand{\bprop}{\begin{prop}}
 \newcommand{\eprop}{\end{prop}}
 \newcommand{\bcor}{\begin{cor}}
 \newcommand{\ecor}{\end{cor}}
 \newcommand{\bnotn}{\begin{notn}}
 \newcommand{\enotn}{\end{notn}}
 \newcommand{\bpb}{\begin{pb}}
 \newcommand{\epb}{\end{pb}}
 \newcommand{\bform}{\begin{form}}
 \newcommand{\eform}{\end{form}}
 \newcommand{\brk}{\begin{rk}}
 \newcommand{\erk}{\end{rk}}
 \newcommand{\bcom}{\begin{com}}
 \newcommand{\ecom}{\end{com}}
 \newcommand{\bex}{\begin{ex}}
 \newcommand{\eex}{\end{ex}}
 \newcommand{\bpf}{\begin{proof}}
 \newcommand{\epf}{\end{proof}}

\newcommand{\cC}{\mathcal{C}}

\newcommand{\cE}{\mathcal{E}}

\newcommand{\cK}{\mathcal{K}}

\newcommand{\cV}{\mathcal{V}}

\newcommand{\bC}{\mathbb{C}}

\newcommand{\bK}{\mathbb{K}}

\newcommand{\bR}{\mathbb{R}}

\newcommand{\be}{\begin{equation}}
\newcommand{\ee}{\end{equation}}
\newcommand{\bal}{\begin{align}}
\newcommand{\eal}{\end{align}}
\newcommand{\ba}{\begin{align*}}
\newcommand{\ea}{\end{align*}}
\newcommand{\bmx}{\begin{matrix}}
\newcommand{\emx}{\end{matrix}}
\newcommand{\bbmx}{\begin{bmatrix}}
\newcommand{\ebmx}{\end{bmatrix}}
\newcommand{\bpmx}{\begin{pmatrix}}
\newcommand{\epmx}{\end{pmatrix}}
\newcommand{\bvmx}{\begin{vmatrix}}
\newcommand{\evmx}{\end{vmatrix}}

\newcommand{\ol}{\overline}
\newcommand{\wh}{\widehat}
\newcommand{\wt}{\widetilde}
\newcommand{\f}{\frac}
\newcommand{\df}{\dfrac}

\newcommand{\setm}{\setminus}

\newcommand{\Id}{\mathrm{Id}}

\newcommand{\argmin}{{\rm argmin}\,}
\newcommand{\argmax}{{\rm argmax}\,}
\newcommand{\minimize}[1]{\underset{#1}{\rm minimize}\,}

\newcommand{\la}{\lambda}
\newcommand{\La}{\Lambda}
\newcommand{\eps}{\varepsilon}
  
\usepackage{natbib}
\title{\vspace{-12mm}Optimal Recovery from Inaccurate Data in Hilbert Spaces:\\
Regularize, but what of the Parameter?

\medskip\hrule height 1.2pt \vspace{-6mm}}
\author{Simon Foucart\footnote{S. F. is supported by grants from the NSF (CCF-1934904, DMS-2053172) and from the ONR (N00014-20-1-2787).} \, and Chunyang Liao  --- Texas A\&M University}
\date{\vspace{-6mm}\rule{100mm}{0.8pt}}

\newcommand\shorttitle{Optimal Recovery from Inaccurate Data in Hilbert Spaces}
\newcommand\authors{S. Foucart, C. Liao}

\begin{document}
\maketitle

%% Add abstract, keywords, and AMS classification
\vspace{-15mm}
\begin{abstract}
In Optimal Recovery,
the task of learning a function from observational data is tackled deterministically by adopting a worst-case perspective tied to
an explicit model assumption made on the functions to be learned.
Working in the framework of Hilbert spaces,
this article considers a model assumption based on approximability.
It also incorporates observational inaccuracies modeled via 
additive errors bounded in $\ell_2$.
Earlier works have demonstrated that regularization provide algorithms that are optimal in this situation,
but did not fully identify the desired hyperparameter.
This article fills the gap in both a local scenario and a global scenario.
In the local scenario,
which amounts to the determination of Chebyshev centers,
the semidefinite recipe of Beck and Eldar
(legitimately valid in the complex setting only)
is complemented by a more direct approach,
with the proviso that the observational functionals have orthonormal representers.
In the said approach,
the desired parameter is the solution to an equation that can be resolved via standard methods.
In the global scenario, 
where linear algorithms rule,
the parameter elusive in the works of Micchelli et al.
is found as the byproduct of a semidefinite program.
Additionally and quite surprisingly,
in case of observational functionals with orthonormal representers,
it is established that any regularization parameter is optimal.
\end{abstract}



\vspace{-5mm}
\begin{center}
\rule{100mm}{0.8pt}
\end{center}

\section{Introduction}

\subsection{Background on Optimal Recovery}

This article is concerned with a central problem in Data Science,
namely:
a function $f$ is acquired through point evaluations
\be
\label{YIdeal}
y_i = f(x^{(i)}),
\qquad i=1,\ldots,m,
\ee
and these data should be used to learn $f$---or to recover it,  with the terminology preferred in this article.
Importantly,
the evaluation points $x^{(1)},\ldots,x^{(m)}$ are considered fixed entities in our scenario:
they cannot be chosen in a favorable way, as in Information-Based Complexity \citep*{novak2008tractability},
nor do they occur as independent realizations of a random variable, as in Statistical Learning Theory \citep*{friedman2001elements}.
In particular,  without an underlying probability distribution, the performance of the recovery process cannot be assessed via generalization error.
Instead, it is assessed via a notion of worst-case error, central to the theory of Optimal Recovery \citep*{MicRiv}.

To outline this theory,  we make the framework slightly more abstract.
Precisely,  given a normed space $F$,
the unknown function is replaced by an element $f \in F$.
This element is accessible only through 
{\em a priori} information expressing an educated belief about $f$
and {\em a posteriori} information akin to \eqref{YIdeal}.
In other words,  our partial knowledge about $f$ is summed up via\vspace{-5mm}
\begin{itemize}
\item the fact that $f \in \cK$ for a subset $\cK$ of $F$ called a model set; 
\item the observational data $y_i = \la_i(f)$, $i = 1,\ldots,m$, for some linear functionals $\la_1,\ldots,\la_m \in F^*$
making up the observation map $\La: g \in F \mapsto [\la_1(g); \ldots; \la_m(g)] \in \bR^m$.
\end{itemize}
We wish to approximate $f$ by some $\wh{f} \in F$ produced using this partial knowledge of $f$.
Since the error $\|f - \wh{f}\|$ involves the unknown $f$,
which is only accessible via $f \in \cK$ and $\La(f)=y$,
we take a worst-case perspective leading to the local worst-case error
\be
\label{LWCE_exact}
{\rm lwce}(y,\wh{f}) := \sup_{\substack{f \in \cK \\ \La(f) = y}} \|f - \wh{f}\|.
\ee
Our objective consists in finding an element $\wh{f}$ that minimizes ${\rm lwce}(y,\wh{f})$.
Such an $\wh{f}$ can be described, almost tautologically,
as a center of a smallest ball containing $\cK \cap \La^{-1}(\{ y \})$.
It is called a Chebyshev center of this set of model- and data-consistent elements.
This remark, however, does not come with any practical construction of a Chebyshev center.

The term local was used above to make a distinction with the global worst-case error of a recovery map $\Delta (=\Delta_\cK): \bR^m \to F$, defined as
\be
\label{GWCE_exact}
{\rm gwce}(\Delta) := \sup_{y \in \La(\cK)} {\rm lwce}(y,\Delta(y)) 
= \sup_{f \in \cK} \|f - \Delta(\La(f))\|.
\ee
The minimal value of ${\rm gwce}(\Delta)$  is called the intrinsic error 
(of the observation map $\La$ over the model set $\cK$)
and the maps $\Delta$ that achieve this minimal value are called globally optimal recovery maps.
Our objective consists in constructing such maps---of course,  the map that assigns to $y$ a Chebyshev center of $\cK \cap \La^{-1}(\{ y \})$ is one of them,
but it may be impractical.
By contrast,
for model sets that are convex and symmetric,
the existence of linear maps among the set of globally optimal recovery maps
is guaranteed by fundamental results from Optimal Recovery in at least two settings:
when $F$ is a Hilbert space 
and when $F$ is an arbitrary normed space but the full recovery of $f$ gives way to the recovery of a quantity of interest $Q(f)$,
$Q$ being a linear functional.
We refer the readers to \citep*[Chapter 9]{BookDS} for details.

\subsection{The specific problem}

The problem solved in this article is a quintessential Optimal Recovery problem---its specificity lies in the particular model set and in the incorporation of errors in the observation process.
The underlying normed  space $F$ is a Hilbert space 
and is therefore denoted by $H$ from now on.
Reproducing kernel Hilbert spaces,
whose usage is widespread in Data Science \citep*{scholkopf2002learning},
 are of particular interest as point evaluations of type \eqref{YIdeal} make perfect sense there.
 
Concerning the model set,
 we concentrate on an approximation-based choice that is increasingly scrutinized, see e.g.  \citep*{maday2015parameterized},
\citep*{devore2017data} and \citep*{cohen2020nonlinear}.
Depending on a linear subspace $\cV$ of $H$
and on a parameter $\eps > 0$,
it takes the form
$$
\cK = \{ f \in H: {\rm dist}(f,\cV) \le \eps \}.
$$
\cite*{binev2017data} completely solved the Optimal Recovery problem with exact data in this situation (locally and globally).
Precisely,  they showed that the solution $\wh{f}$ to
\be
\label{OptProgBinev}
\minimize{f \in H} \; {\rm dist}(f,\cV)
\qquad \mbox{s.to} \quad \La(f)=y,
\ee
which clearly  belongs to the model- and data-consistent set $\cK \cap \La^{-1}(\{y\})$, 
turns out to be its Chebyshev center.
Moreover,  with $P_\cV$ and $P_{\cV^\perp}$ denoting the orthogonal projectors onto $\cV$ and onto the orthogonal complement $\cV^\perp$ of $\cV$,  the fact that ${\rm dist}(f,\cV) = \|f - P_{\cV}f\| = \|P_{\cV^\perp} f\|$ makes 
the optimization program \eqref{OptProgBinev} tractable.
It can actually be seen that $\Delta: y \mapsto \wh{f}$ is a linear map.
This is a significant advantage
because $\Delta$ can then be precomputed in an offline stage knowing only $\cV$ and $\La$
 and  the program~\eqref{OptProgBinev} need not be solved afresh for each new data $y \in \bR^m$ arriving in an online stage.

Concerning the observation process,
instead of exact data $y = \La(f) \in \bR^m$,
it is now assumed that
$$
y = \La(f) + e \in \bR^m
$$
for some unknown error vector $e \in \bR^m$.
This error vector is not modeled as random noise but through the deterministic $\ell_2$-bound $\|e\|_2 \le \eta$.
Although other $\ell_p$-norms can the considered for the optimal recovery of $Q(f_0)$ when $Q$ is a linear functional on an arbitrary normed space $F$ 
(see \citep*{ettehad2020instances}),
here the arguments rely critically on $\bR^m$ being endowed with the $\ell_2$-norm.
It will be written simply as $\|\cdot\|$ below, hoping that it does not create confusion with the Hilbert norm on $H$.

For our specific problem,  the worst-case recovery errors \eqref{LWCE_exact} and \eqref{GWCE_exact} need to be adjusted.
The local worst-case recovery error at $y$ for $\wh{f}$ becomes
$$
{\rm lwce}(y,\wh{f}) = \sup_{\substack{ \|P_{\cV^\perp} f\| \le \eps \\ \|\La(f)-y\| \le \eta }}
\|f-\wh{f}\|.
$$
As for the global worst-case error of $\Delta: \bR^m \to H$,
it reads
$$
{\rm gwce}(\Delta) = \sup_{\substack{ \|P_{\cV^\perp} f\| \le \eps \\ \|e\| \le \eta }}
\| f - \Delta(\La(f)+e) \|.
$$
Note that both worst-case errors are infinite if one can find a nonzero $h$  in $\cV \cap \ker(\La)$.
Indeed, the element $f_t := f + t h$, $t \in \bR$,
obeys $\|P_{\cV^\perp} f_t \| = \|P_{\cV^\perp} f \|  \le \eps$ and $\|y-\La(f_t)\| = \|y-\La(f)\| \le \eta$,
so for instance ${\rm lwce}(y,\wh{f}) \ge \sup_{t \in \bR} \|f_t -\wh{f} \| = +\infty$.
Thus, we always make the assumption that 
\be
\label{Assum_VK}
\cV \cap \ker(\La) = \{0\}.
\ee
We keep in mind that the latter forces $n:= \dim(\cV) \le m$,
as can be seen by dimension arguments.
With $\La^*$ denoting the Hermitian adjoint of $\La$,
another assumption that we sometimes make reads
\be
\label{Assum_TF}
\La \La^* = \Id_{\bR^m}.
\ee
This is not extremely stringent:
assuming the surjectivity of $\La$ is quite natural,
otherwise certain observations need not be collected;
then the map $\La$ can be preprocessed into another map $\wt{\La}$ satisfying $\wt{\La} \wt{\La}^* = \Id_{\bR^m}$ by setting $\wt{\La} = (\La \La^*)^{-1/2} \La$.
Incidentally, if $u_1,\ldots,u_m \in H$ represent the Riesz representers of the observation functionals $\la_1,\ldots,\la_m \in H^*$,
characterized by $\langle u_i, f \rangle  = \la_i(f)$ for all $f \in H$,
then the assumption \eqref{Assum_TF} is equivalent to the orthonormality of the system $(u_1,\ldots,u_m)$.
In a reproducing kernel Hilbert space with kernel $K$,
if the $\la_i$'s are point evaluations at some $x^{(i)}$'s,
so that $u_i = K(\cdot,x^{(i)})$,
then \eqref{Assum_TF} is equivalent to $K(x^{(i)},x^{(j)}) = \delta_{i,j}$ for all $i,j =1,\ldots,m$.
This occurs e.g. for the Paley--Wiener space of functions with Fourier transform supported on $[-\pi,\pi]$
when the evaluations points come from an integer grid,
since the kernel is given by $K(x,x') = {\rm sinc}(\pi(x-x'))$,
$x,x' \in \bR$.

\subsection{Main results}

There are previous works on Optimal Recovery in Hilbert spaces in the presence of observation error bounded in $\ell_2$.
Notably, 
%\citep*{beck2006strong} and 
\citep*{beck2007regularization} dealt with the local setting,
while \citep*{micchelli1979inaccurate} and \citep*{micchelli1993optimal} dealt with the global setting.
These works underline the importance of regularization,
which is prominent in many other settings  \citep*{chen2002different}.
They establish that the optimal recovery maps are obtained by solving the unconstrained program
\be
\label{Reg}
\minimize{f \in H} \; (1-\tau) \|P_{\cV^\perp} f\|^2 + \tau \| \La f - y \|^2
\ee
for some $\tau \in [0,1]$.
It is the precise choice of this regularization parameter $\tau$ which is the purpose of this article.
Assuming from now on that $H$ is finite dimensional\footnote{It is likely that the results are still valid in the infinite-dimensional case.
But then it is unclear how to solve \eqref{EqEig} and \eqref{SDP4Tau} numerically, so the infinite-dimensional case is not given proper scrutiny in the rest of the article.},
we provide a  complete (almost) picture of the local and global Optimal Recovery solutions,
as summarized in the four points below, three of them being new:
\vspace{-3mm}
\begin{enumerate}

\item[{\bf L1.}] With $H$ restricted here to be a complex Hilbert space,
the Chebyshev center of the set $\{ f \in H: \|P_{\cV^\perp} f \| \le \eps, \| \La f - y \| \le \eta \}$ is the minimizer of \eqref{Reg} for the choice $\tau = d_\sharp/(c_\sharp + d_\sharp )$,
where $c_\sharp,d_\sharp$ are solutions to the semidefinite program\footnote{In the statement of this semidefinite program and elsewhere,
the notation $T \succeq 0$ means that an operator $T$ is positive semidefinite on $H$,
i.e., that $\langle Tf,f \rangle \ge 0$ for all $f \in H$.}
\begin{align*}
\minimize{c,d,t \ge 0} \; \eps^2 c + (\eta^2-\|y\|^2) d + t
& \qquad \mbox{s.to} & & 
c P_{\cV^\perp} + d \La^* \La \succeq \Id,\\
&  \qquad \mbox{and} & & \bbmx
c P_{\cV^\perp} + d \La^* \La & | & -d \La^* y\\
\hline
-d(\La^* y)^*  & | & t
\ebmx \succeq 0.\\
\end{align*}

\item[{\bf L2.}] Under the orthonormal observations assumption \eqref{Assum_TF} but without the above restriction on $H$,
the Chebyshev center of the set $\{ f \in H: \|P_{\cV^\perp} f \| \le \eps, \| \La f - y \| \le \eta \}$ is the minimizer of \eqref{Reg} for the choice $\tau$ that satisfies
\be
\label{EqEig} 
\la_{\min} ((1-\tau)P_{\cV^\perp} + \tau \La^* \La)
= \f{(1-\tau)^2 \eps^2 - \tau^2 \eta^2}{(1-\tau) \eps^2 - \tau \eta^2 + (1-\tau) \tau (1-2 \tau)  \delta^2},
\ee
where $\delta$ is precomputed as $\delta = \min\{ \|P_{\cV^\perp}f\|: \La f = y \} = \min \{ \|\La f - y\|: f \in \cV \}$.
For the distinct case $\cV = \{0\}$, the best choice of parameter is more simply $\tau = \max\{1 - \eta / \|y\|,0\}$.

\item[{\bf G1.}] A globally optimal recovery map is provided by the linear map sending $y \in \bR^m$ to the minimizer of \eqref{Reg} with parameter $\tau = d_\flat/(c_\flat+d_\flat)$, 
where $c_\flat,d_\flat$ are solutions to the semidefinite program
\be
\label{SDP4Tau}
\minimize{c,d} \; \eps^2 c + \eta^2 d
\qquad \mbox{s.to}  \quad
c P_{\cV^\perp} + d \La^* \La \succeq \Id.
\ee

\item[{\bf G2.}] Under the orthonormal observations assumption \eqref{Assum_TF},
the linear map sending $y \in \bR^m$ to the minimizer of \eqref{Reg} 
is a globally optimal recovery map for any choice of parameter $\tau \in [0,1]$.

\end{enumerate}

Before entering the technicalities,
a few comments are in order to put these results in context.
Item~{\bf L1} is the result of \citep*{beck2007regularization}
(see Corollary 3.2 there)
adapted to our situation.
It relies on an extension of the S-lemma involving two quadratic constraints.
This extension is valid in the complex finite-dimensional setting,
but not necessarily in the real setting, 
hence the restriction on $H$
(this does not preclude the validity of the result in the real setting, though).
It is worth pointing out the nonlinearity of the map that sends $y \in \bR^m$ to the above Chebyshev center.
Incidentally,  we can safely talk about the Chebyshev center,
because it is known \citep*{garkavi1962optimal} that a bounded set in a uniformly convex Banach space has exactly one Chebyshev center.
A sketch of the argument adapted to our situation is presented in the appendix.
\label{PageCCUnique}

For item~{\bf L2},
working with an observation map $\La$ satisfying $\La \La^* = \Id_{\bR^m}$
allows us to construct the Chebyshev center even in the setting of a real Hilbert space.
This is possible because our argument does not rely on the extension of the S-lemma---it just uses the obvious implication. 
As for equation \eqref{EqEig},
it is easily solved using the bisection method or the Newton/secant method.
Moreover, it gives some insight on the value of the optimal parameter $\tau$.
For instance, the proof reveals that $\tau$ is always between $1/2$ and $\eps/(\eps + \eta)$.
When $\eps \ge \eta$, say, the optimal parameter should then satisfy $\tau \ge 1/2$,
which is somewhat intuitive:
$\eps \ge \eta$ means that there is more model mismatch than data mismatch,
so the regularization should penalize model fidelity less than data fidelity
by taking $1-\tau \le \tau$, i.e., $\tau \ge 1/2$.
As an aside,  we point out that, here too, the map that sends $y \in \bR^m$ to the Chebyshev center is not a linear map---if it was,  then the optimal parameter should be independent of $y$.

In contrast,
the globally optimal recovery map of item~{\bf G1} is linear.
It is one of several globally optimal recovery maps,
since the locally optimal one (which is nonlinear) is also globally optimal.
However, as revealed in the reproducible\footnote{{\sc matlab} and Python files illustrating the findings of this article are located at \url{https://github.com/foucart/COR}.} accompanying this article,
it is in general the only regularization map that turns out to be globally optimal. 
The fact that regularization produces globally optimal recovery maps was recognized by Micchelli,
who wrote in the abstract of \citep*{micchelli1993optimal}
that ``{\em the regularization parameter must be chosen with care}''.
However, a recipe for selecting the parameter was not given there,
except on a specific example.
The closest to a nonexhaustive search is found in \citep*[Lemma 2.6.2]{plaskota1996noisy} for the case $\cV = \{0\}$,
but even this result does not translate into a numerically tractable recipe.
The selection stemming from \eqref{SDP4Tau} does,
at least when $H$ is finite-dimensional, which is assumed here.
Semidefinite programs can indeed be solved
in {\sc matlab} using {\sf CVX} \citep*{CVX} and in Python using {\sf CVXPY} \citep*{CVXPY}.
%The argument leading to \eqref{SDP4Tau} is also slicker.

Finally, a surprise arises in item~{\bf G2}.
Working with an observation map $\La$ satisfying $\La \La^* = \Id_{\bR^m}$,
the latter indeed reveals that {\em the regularization parameter does not need to be chosen with care} after all,
since regularization maps are globally optimal
no matter how the parameter $\tau \in [0,1]$ is chosen.
The precise interpretation of the choices $\tau = 0$ and $\tau=1$
will be elucidated later.

The rest of this article is organized as follows.
Section~\ref{SecTech} gathers some auxiliary results that are used in the proofs of the main results.
Section~\ref{SecLoc} elucidates item~{\bf L1} and establishes item~{\bf L2}---in other words, it is concerned with local optimality.
Section~\ref{SecGlob}, which is concerned with global optimality,
is the place where items~{\bf G1} and {\bf G2} are proved.
Lastly, a short appendix containing some side information is included after the bibliography.

\section{Technical Preparation}
\label{SecTech}

This section establishes (or recalls) a few results that we isolate here in order not to disrupt the flow of subsequent arguments. 

\subsection{S-lemma and S-procedure}

Loosely speaking, the S-procedure is a relaxation technique expressing the fact that a quadratic inequality is a consequence of some quadratic constraints.
In case of a single quadratic constraint,
the relaxation turns out to be exact.
This result, known as the S-lemma, can be stated as follows:
given quadratic functions $q_0$ and $q_1$ defined on $\bK^N$, with $\bK = \bR$ or $\bK = \bC$,
$$
[q_0(x) \le 0
\; \mbox{ whenever } q_1(x) \le 0]
\iff [\mbox{there exists } a \ge 0: q_0 \le a q_1],
$$
provided $q_1(\wt{x})<0$ for some $\wt{x} \in \bK^N$.
With more than one quadratic constraint,
 $q_1,\ldots,q_k$, say,
$q_0(x) \le 0$ whenever $q_1(x)\le 0,\ldots,q_k(x)\le 0$
is still a consequence of $q_0 \le a_1 q_1 + \cdots + a_k q_k$ 
for some $a_1,\ldots,a_k \ge 0$,
but the reverse implication does not hold anymore.
There is a subtlety when $k=2$,
as the reverse implication holds for $\bK = \bC$
but not for $\bK = \bR$,  see \citep[Section~3]{polik2007survey}.
However, if the quadratic constraints do not feature linear terms, 
then the reverse implication holds for $k=2$ also when $\bK = \bR$.
Since this result of \citep[Theorem 4.1]{polyak1998convexity} is to be invoked later, 
we state it formally below.
 
\bthm
\label{ThmPolyak}
Suppose that $N \ge 3$ and that quadratic functions $q_0,q_1,q_2$ on $\bR^N$ 
take the form $q_i(x) = \langle A_i x,x \rangle + \alpha_i$ for symmetric matrices $A_0,A_1,A_2 \in \bR^{N \times N}$ and scalars $\alpha_0,\alpha_1,\alpha_2 \in \bR$.
Then
$$
[q_0(x) \le 0
\; \mbox{ whenever } q_1(x) \le 0 \mbox{ and } q_2(x) \le 0]
\iff 
[\mbox{there exist }a_1 ,  a_2 \ge 0:
q_0 \le a_1 q_1 + a_2 q_2],
$$
provided $q_1(\wt{x})<0$ and $q_2(\wt{x})<0$ for some $\wt{x} \in \bR^N$
and $b_1 A_1 + b_2 A_2 \succ 0$ for some $b_1,b_2 \in \bR$.
\ethm

\subsection{Regularization}

In this subsection,
we take a closer look at the regularization program \eqref{Reg}.
The result below shows that its solution depends linearly on $y \in \bR^m$.
In fact, the result covers a slightly more general program
and the linearity claim follows by taking $R = P_{\cV^\perp}$, $r=0$,  $S = \La$, and $s=y$.

\bprop
\label{PropGenExpr}
Let $R,S$ be linear maps from $H$ into other Hilbert spaces containing points $r,s$, respectively.
For $\tau \in (0,1)$,
the optimization program
\be
\label{GenReg}
\minimize{f \in H} \; (1-\tau) \|Rf-r\|^2 + \tau \|Sf-s\|^2
\ee
has solutions $f_\tau \in H $ characterized by 
\be
\label{CharaGenReg}
\big( (1-\tau)R^* R+ \tau S^*S \big) f_\tau = (1-\tau)R^*r + \tau S^*s.
\ee
Moreover, if $\ker(R) \cap \ker(S) = \{0\}$, then $f_\tau$ is uniquely given by
\be
\label{ExprGenRegInv}
 f_\tau = \big( (1-\tau)R^*R + \tau S^*S \big)^{-1} \big( (1-\tau)R^*r + \tau S^*s \big).
\ee
\eprop

\bpf
The program \eqref{GenReg} can be interpreted as a standard least squares problem, namely as
$$
\minimize{f \in H} \left\| \bbmx \sqrt{1-\tau}R \\ \hline \sqrt{\tau} S \ebmx f - \bbmx \sqrt{1-\tau} r\\ \hline \sqrt{\tau} s \ebmx \right\|^2.
$$
According to the normal equations, 
its solutions $f_\tau$ are characterized by  
$$
\bbmx \sqrt{1-\tau}R \\ \hline \sqrt{\tau} S \ebmx^* \bbmx \sqrt{1-\tau}R \\ \hline \sqrt{\tau} S \ebmx f_\tau
= \bbmx \sqrt{1-\tau}R \\ \hline \sqrt{\tau} S \ebmx^* \bbmx \sqrt{1-\tau} r\\ \hline \sqrt{\tau} s \ebmx,
$$
which is a rewriting of \eqref{CharaGenReg}.
Next, if $\ker(R) \cap \ker(S) = \{0\}$, then
$$
\langle  ((1-\tau) R^* R + \tau S^* S) f, f \rangle
= (1-\tau) \|Rf\|^2 + \tau \|Sf\|^2 \ge 0,
$$
with equality only possible when $f \in \ker(R) \cap \ker(S)$,
i.e.,  $f=0$.
This shows that $(1-\tau) R^* R + \tau S^* S$ is positive definite,
and hence invertible,
which allows us to write \eqref{ExprGenRegInv} as a consequence of \eqref{CharaGenReg}.
\epf

The expression \eqref{ExprGenRegInv} is not always the most convenient one.
Under extra conditions on $R$ and~$S$,
we shall see that $f_\tau$, $\tau \in [0,1]$, can in fact be expressed as the convex combination $f_\tau = (1-\tau) f_0 + \tau f_1$.
The elements $f_0$ and $f_1$ should be interpreted\footnote{Intuitively,  
the solution to the program \eqref{GenReg} written as the minimization of $\|Rf-r\|^2 + (\tau/(1-\tau))\|Sf-s\|^2$ becomes,  as $\tau \to 1$,  the mininizer of $\|Rf-r\|^2$ subject to $\|Sf-s\|^2=0$.
This explains the interpretation of $f_1$.
A~similar argument explains the interpretation of $f_0$.} as 
\begin{align*}
f_0  = \underset{f \in H}{\argmin} \|Sf-s\| & \qquad  \mbox{s.to } \; Rf=r,\\
f_1  = \underset{f \in H}{\argmin} \|Rf-r\| & \qquad \mbox{s.to } \; Sf=s.\\
\end{align*} 
The requirements that $r \in {\rm ran}(R)$ and $s \in {\rm ran}(S)$ need to be imposed for $f_0$ and $f_1$ to even exist
and the condition $\ker(R) \cap \ker(S) = \{0\}$ easily guarantees that $f_0$ and $f_1$ are unique.
They obey
\be
\label{IdenF0F1}
R f_0 = r,
\quad S^*(S f_0 - s) \in \ker(R)^\perp ,
\quad S f_1 = s,
\quad R^*(R f_1 - r) \in \ker(S)^\perp.
\ee
For instance, the identity $R f_0 = r$ reflects the constraint in the optimization program defining $f_0$,
while $S^*(S f_0 - s) \in \ker(R)^\perp$ is obtained by expanding 
$\|S(f_0 + tu) - s\|^2 \ge \|S(f_0 ) - s\|^2$ around $t=0$ for any $u \in \ker(R)$.
At this point, we are ready to establish our claim under extra conditions on $R$ and $S$,
namely that they are orthogonal projections.
These conditions will be in place when the observation map satisfies $\La \La^* = \Id_{\bR^m}$.
Indeed,  in view of $\|w\|^2 = \langle w, \La \La^* w \rangle =  \|\La^* w\|^2$ for any $w \in \bR^m$,
the regularization program \eqref{Reg} also reads
$$
\minimize{f \in H} \; (1-\tau) \|P_{\cV^\perp} f\|^2 + \tau \| \La^*\La f - \La^* y \|^2,
$$
where both $P_{\cV^\perp}$ and $\La^* \La$ are orthogonal projections.
The result below will then be applied with  $R = P_{\cV^\perp}$, $r=0$,  $S = \La^* \La$, and $s=\La^* y$.

\bprop
\label{PropLinearExpr}
Let $R,S$ be two orthogonal projectors on $H$ such that $\ker(R) \cap \ker(S) = \{0\}$
and let $r\in {\rm ran}(R)$, $s \in {\rm ran}(S)$.
For $\tau \in [0,1]$,
%Assume that $\ker(R)$ and ${\rm ran}(S)$ are finite-dimensional,
%with $(v_1,\ldots,v_k)$ being a basis for $\ker(R)$,
%$(w_1,\ldots,w_\ell)$ being an orthonormal basis for ${\rm ran}(S)$,
%and $C \in \bR^{k \times \ell}$ denoting the cross-gramian with entries $C_{i,j} = \langle v_i, w_j \rangle$.
the solution $f_\tau $ to the optimization program
\be
\label{GenReg2}
\minimize{f \in H} \; (1-\tau) \|Rf-r\|^2 + \tau \|Sf-s\|^2
\ee
satisfies
\be
\label{ExprGenRegAff}
f_\tau = 
(1-\tau) f_0 + \tau f_1.
%\sum_{i=1}^k a_i v_i + \tau \sum_{j=1}^\ell b_j w_j + r.
\ee
%The coefficient vectors $a \in \bR^k$ and $b \in \bR^\ell$  depend linearly on $r$ and $s$ according to
%\be
%\label{Defab}
%a = \big( C C^\top )^{-1} C [s-r]_w
%\qquad \mbox{and} \qquad
%b = [s-r]_w - C^\top a,
%\ee
%where $[s-r]_w \in \bR^\ell$ has entries $\langle s-r,w_j \rangle$ for $j=1,\ldots,\ell$. 
Moreover, one has
\be
\label{NormRFtau-r}
\|Rf_\tau - r\| = \tau \|f_1-f_0\|
\qquad \mbox{and} \qquad
\|Sf_\tau-s\|=(1-\tau)\|f_1-f_0\|.
\ee
\eprop

\bpf
Taking the extra conditions on $R$ and $S$ into account,
the identities  \eqref{IdenF0F1} read
$$
R f_0 = r,
\quad S f_0 - s \in {\rm ran} (R),
\quad S f_1 = s,
\quad R f_1 - r \in {\rm ran} (S).
$$
In this form, they imply that 
\begin{align}
\nonumber
\langle S(f_0-f_1),(R-S)(f_0-f_1) \rangle
& = \langle Sf_0-s,R(f_0-f_1) \rangle - \langle Sf_0-s,S(f_0-f_1) \rangle\\
\nonumber
& =  \langle Sf_0-s,f_0-f_1 \rangle - \langle Sf_0-s,f_0-f_1 \rangle \\
\label{SRF0F1}
& = 0.
\end{align}
In a similar fashion, by exchanging the roles of $R$ and $S$,
and consequently also of $f_0$ and $f_1$,
we have
$\langle R(f_1-f_0),(S-R)(f_1-f_0) \rangle = 0$,
i.e.,  $\langle R(f_0-f_1),(R-S)(f_0-f_1) \rangle = 0$.
Subtracting \eqref{SRF0F1} from the latter yields $\|(R-S)(f_0-f_1)\|^2 = 0$,
in other words $R(f_0-f_1) = S(f_0-f_1)$.
Then, 
the element $h:= f_0-f_1 - R(f_0-f_1) = f_0-f_1 - S(f_0-f_1)$ belongs to $\ker(R) \cap \ker(S)$, so  that $h=0$.
In summary, we have established that
\be
\label{f0Nf1}
R(f_0- f_1) = S(f_0-f_1) = f_0-f_1.
\ee
From here,  we can deduce the two parts of the proposition.
For the first part, we notice that
\begin{align*}
\big( (1-\tau) R + \tau S \big) \big(  (1-\tau) f_0 + \tau f_1  \big)
& = (1-\tau)^2 Rf_0 + (1-\tau) \tau (Rf_1 + S f_0) + \tau^2 Sf_1\\
& = (1-\tau)^2 Rf_0 + (1-\tau) \tau (Rf_0 + S f_1) + \tau^2 Sf_1\\
& = (1-\tau) Rf_0 + \tau S f_1 \\
& = (1-\tau) r + \tau s, 
\end{align*}
which shows that
$(1-\tau)f_0 + \tau f_1$ satisfies the relation \eqref{CharaGenReg} characterizing the minimizer $f_\tau$ of \eqref{GenReg2},
so that $f_\tau = (1-\tau)f_0 + \tau f_1$,
as announced in \eqref{ExprGenRegAff}.
For  the second part, 
we notice that
$$
Rf_\tau - r = (1-\tau) R f_0 + \tau Rf_1 - Rf_0 = \tau R(f_1-f_0) = \tau( f_1-f_0),
$$
so the first equality of \eqref{NormRFtau-r} follows by taking the norm.
The second equality of \eqref{NormRFtau-r} is derived in a similar fashion.
\epf

We complement Proposition \ref{PropLinearExpr} with a few additional pieces of information.

\brk
Under the assumptions of Proposition \ref{PropLinearExpr},
the solution $f_\tau$ to \eqref{GenReg2} is also solution to
$$
\minimize{f \in H} \; \max\big\{  (1-\tau) \|Rf-r\|, \tau \|Sf-s\|  \big\}.
$$
Indeed, at $f=f_\tau$,
the squared objective function equals $(1-\tau)^2 \tau^2 \|f_1-f_0\|^2$,
while at an arbitrary $f \in H$,
it satisfies
\begin{align*}
 \max\big\{  (1-\tau)^2 \|Rf-r\|^2, \tau^2 \|Sf-s\|^2  \big\}
 & \ge \f{1}{\df{1}{1-\tau} + \df{1}{\tau}} \big( (1-\tau) \|Rf-r\|^2 + \tau \|Sf-s\|^2  \big)\\
&  \ge (1-\tau) \tau \big( (1-\tau) \|Rf_\tau-r\|^2 + \tau \|Sf_\tau-s\|^2  \big)\\
& = (1-\tau) \tau \big( (1-\tau) \tau^2 \|f_1-f_0\|^2 + \tau (1-\tau)^2 \|f_1-f_0\|^2  \big)\\
& = (1-\tau)^2 \tau^2 \|f_1-f_0\|^2.
\end{align*}
In the case $R = P_{\cV^\perp}$, $r=0$,  $S = \La^* \La$, and $s=\La^* y$,
the choice $\tau = \eps/(\eps + \eta)$ is quite relevant,
since the above optimization program becomes equivalent to
$$
\minimize{f \in H} \; \max\bigg\{  \f{1}{\eps} \|P_{\cV^\perp}f\|,   \f{1}{\eta} \|\La f-y\|  \bigg\}.
$$
Its solution is clearly in the model- and data-consistent set $\{ f \in H: \|P_{\cV^\perp}f\| \le \eps,  \|\La f-y\| \le \eta \}$.
In~fact,  this could have been a natural guess for its Chebyshev center,
but item~{\bf L2} reveals the invalidity of such a guess.
Nonetheless, the special parameter $\tau = \eps/(\eps + \eta)$ will make a reappearance in the argument leading to item~{\bf L2}.
\erk

\brk
The proof of Proposition \ref{PropLinearExpr} showcased the important identities
$Rf_0 = r$, $S f_1=s$, and $R(f_0-f_1) = S(f_0-f_1)=f_0-f_1$.
In the case $R = P_{\cV^\perp}$, $r=0$,  $S = \La^* \La$, and $s=\La^* y$,
if $\Delta_\tau$ denotes the recovery map assigning to $y \in \bR^m$ the solution $f_\tau$ to the regularization program \eqref{Reg}, these identities read,
when $\La \La^* = \Id_{\bR^m}$,
\be
\label{PropertiesDeltaUnderTF}
P_{\cV^\perp} \Delta_0 =0,
\qquad \La^* \La \Delta_1 =  \La^*,
\qquad  P_{\cV^\perp}( \Delta_0 - \Delta_1) = 
\La^* \La(\Delta_0 - \Delta_1) = \Delta_0 - \Delta_1.
\ee
\erk

\brk
Considering again the case $R = P_{\cV^\perp}$, $r=0$,  $S = \La^* \La$, and $s=\La^* y$,
Proposition \ref{PropLinearExpr} implies that $f_\tau \in \cV + {\rm ran}(\La^*)$ for any $\tau \in [0,1]$,
given that the latter holds for $\tau=0$ and for $\tau=1$.
For $\tau=0$,
this is because the constraint $P_{\cV^\perp} f = 0$ of the optimization program defining $f_0$ imposes $f_0 \in  \cV$.
For $\tau=1$,
this is a result established e.g.  in \citep*[Theorem 2]{foucart2020learning}.
The said result also provides an efficient way to compute the solution $f_\tau$ of \eqref{Reg} even when $H$ is infinite dimensional, as stated in the appendix.
\label{PageCompReg}
\erk

\section{Local Optimality}
\label{SecLoc}

Our goal in this section is to determine locally optimal recovery maps.
In other words,  the section is concerned with Chebyshev centers.
We start by considering the situation of an arbitrary observation map $\La$,
but with a restriction on the space $H$.
Next, lifting this restriction on $H$,
we refine the result in the particular case of an observation map satisfying $\La \La^* = \Id_{\bR^m}$.

\subsection{Arbitrary observations}

In this subsection, we reproduce a result from \citep{beck2007regularization},
albeit with different notation,
and explain how it implies the statement of item~{\bf L1}.
The result in question, namely Corollary~3.2, relies on the S-procedure with two constraints,
and as such cannot be claimed in the real setting.

\bthm
\label{ThmLocGenObs}
Let $H$ be a complex Hilbert space.
Let $R,S$ be two linear maps from $H$ into other Hilbert spaces containing points $r,s$, respectively.
Suppose the existence of $\wt{f} \in H$ such that $\|R \wt{f} - r \| < \eps$ and $\|S \wt{f} - s \| < \eta$
and the existence of $\tau \in [0,1]$ such that $(1-\tau) R^* R + \tau S^* S$ is positive definite.
Then the Chebyshev center of $\{ f \in H: \|R f - r \| \le \eps , \|S f - s \| \le \eta \}$
equals $f_\sharp = \big( c_\sharp R^*R + d_\sharp S^*S \big)^{-1} (c_\sharp R^*r+d_\sharp S^* s)$,
where $c_\sharp,d_\sharp$ are solutions to
\begin{align*}
\minimize{c,d,t \ge 0} \; (\eps^2-\|r\|^2) c + (\eta^2-\|s\|^2) d + t
& \qquad \mbox{s.to} & & 
c R^* R + d S^* S \succeq \Id,\\
&  \qquad \mbox{and} & & \bbmx
c R^* R + d S^* S & | & -c R^* r- d  S^* s\\
\hline
-c (R^* r)^* -d(S^* s)^*  & | & t
\ebmx \succeq 0.\\
\end{align*}
\ethm

The statement made in item~{\bf L1} is of course derived by taking $R=P_{\cV^\perp}$, $r=0$, $S = \La$, and $s=y$.
Theorem~\ref{ThmLocGenObs} is indeed applicable, as
$\wt{f}=(f_0+f_1)/2$ satisfies the strict feasibility conditions,
while the positive definiteness condition is not only fulfilled for some $\tau \in [0,1]$, but for all $\tau \in (0,1)$,
since $\langle ( (1-\tau) P_{\cV^\perp} + \tau \La^* \La ) f, f \rangle = (1-\tau) \|P_{\cV^\perp} f\|^2 + \tau \|\La f\|^2 \ge 0$,
with equality only possible if $f \in \cV \cap \ker(\La)$,
i.e., if $f=0$ thanks to the assumption \eqref{Assum_VK}.
We also note that, by virtue of~\eqref{ExprGenRegInv}, the element $f_\sharp$ defined above is nothing else than the regularized solution with parameter $\tau = d_\sharp / (c_\sharp + d_\sharp)$.

\subsection{Orthonormal observations}

In this subsection,
we place ourselves in the situation of an observation map satisfying
$\La \La^* = \Id_{\bR^m}$
and we provide a proof of the statements made item~{\bf L2}.
In fact, we prove some slightly more general results
and {\bf L2} follows by taking 
$R = P_{\cV^\perp}$, $r=0$,  $S = \La^* \La$, and $s= \La^* y$.
Note that we must separate the cases where $R = \Id$ (corresponding to $\cV = \{0\}$)
and where $R$ is a proper orthogonal projection (corresponding to $\cV \not= \{0\}$).
We emphasize that, in each of these two cases,
 the optimal parameter $\tau_\sharp$ is not independent of $y$.
 Therefore,  in view of \eqref{ExprGenRegAff} and of the linear dependence of $f_0$ and $f_1$ on $y$,
 the regularized solution $f_{\tau_\sharp}$ does not depend linearly on $y$.
In other words, the locally optimal recovery map is not a linear map.
The following two simple lemmas will be used to deal with both cases.
%It is worth pointing out that Lemma \ref{LemCheCen} cannot be invoked for arbitrary observations.
%Indeed, the Chebyshev center described in Theorem \ref{ThmLocGenObs} above fails to fulfill the orthogonality conditions \eqref{CondOrtho}.\footnote{\TODO{as shown in the reproducible.}}

\blem
\label{LemCheCen}
Let $R,S$ be two linear maps from $H$ into other Hilbert spaces containing points $r,s$, respectively.
Given $f_\sharp \in H$, let 
$$
h_\sharp \in \underset{h \in H}\argmax \|h\|
\qquad \mbox{s.to } 
\left\{ \bmx 
\|Rf_\sharp - r + Rh\| \hspace{-3mm}& \le \eps ,\\ \|Sf_\sharp -s +Sh\| \hspace{-3mm} & \le \eta.
\emx
\right.
$$
If the orthogonality conditions
\be
\label{CondOrtho}
\langle R^*(Rf_\sharp -r ), h_\sharp \rangle = 0  
\qquad \mbox{and} \qquad  
\langle S^*(Sf_\sharp -s ), h_\sharp \rangle = 0
\ee
are fulfilled,
then $f_\sharp$ is the Chebyshev center of the set $\{ f \in H: \|Rf-r\|\le \eps, \|Sf-s\| \le \eta \},$ 
i.e., for any $g \in H$,
\be
\label{IsCheCenter}
\sup_{\substack{ \|Rf-r\| \le \eps \\ \|Sf-s\| \le \eta }} \|f-g\|
\ge \sup_{\substack{ \|Rf-r\| \le \eps \\ \|Sf-s\| \le \eta }} \|f-f_\sharp\|.
\ee
\elem

\bpf
First, writing $f = f_\sharp + h$,
we easily see that the right-hand side of \eqref{IsCheCenter} reduces to $\|h_\sharp\|$.
Second, let us remark that the orthogonality conditions guarantee that $f_\pm := f_\sharp \pm h_\sharp$ both satisfy 
$\|Rf_\pm-r\| \le \eps$ and $\|Sf_\pm-s\| \le \eta$.
For instance, we have
\be
\|Rf_\pm-r\|^2 = \|Rf_\sharp - r \pm Rh_\sharp\|^2 = \|Rf_\sharp - r\|^2 + \|Rh_\sharp\|^2
= \|Rf_\sharp - r + Rh_\sharp\|^2 \le \eps^2,
\ee
where the latter inequality reflects the feasibility of $h_\sharp$.
Therefore, the left-hand side of \eqref{IsCheCenter} is bounded below by
\be
\max_\pm \|f_\pm - g\| \ge \f{1}{2} \big( \|f_+ -g\| + \|f_- -g\| \big)
\ge \f{1}{2} \|(f_+ -g) - (f_- -g)\| = \f{1}{2} \|2h_\sharp\| = \|h_\sharp\|,
\ee
i.e., by the right-hand side of \eqref{IsCheCenter}. 
\epf

The next lemma somehow relates to the S-procedure.
However, it does not involve the coveted (and usually invalid) equivalence,
but only the straightforward implication.

\blem
\label{LemMaximizer}
Let $R,S$ be two linear maps from $H$ into other Hilbert spaces containing points $r,s$, respectively.
Given $f_\sharp \in H$ and $h_\sharp \in H$, suppose that
\be
\label{CondNorm}
\|Rf_\sharp -r + R h_\sharp\|^2 = \eps^2
\qquad \mbox{and} \qquad
\|Sf_\sharp -s + S h_\sharp\|^2 = \eta^2,
\ee
and that there exist $a,b \ge 0$ such that
\be
\label{Condab1}
a R^* R + b S^* S \succeq \Id
\ee
as well as
\be
\label{Condab2}
a R^*(R f_\sharp - r) + b S^*(Sf_\sharp - s) + (a R^* R + b S^* S ) h_\sharp = h_\sharp.
\ee
Then, one has 
\be
\label{DefHsharp}
h_\sharp \in \underset{h \in H}\argmax \|h\|
\qquad \mbox{s.to } 
\left\{ \bmx 
\|Rf_\sharp - r + Rh\| \hspace{-3mm} &\le \eps ,  \\ \|Sf_\sharp -s +Sh\| \hspace{-3mm} & \le \eta.
\emx
\right.
\ee
\elem

\bpf
By writing the variable in the optimization  program \eqref{DefHsharp} as $h = h_\sharp + g$,
the constraints on~$h$ transform into constraints on~$g$.
Thanks to \eqref{CondNorm}, the latter constraints read
$$
\langle R^* R g,g \rangle + 2 \langle R^*(Rf_\sharp - r + Rh_\sharp), g \rangle \le 0
\quad \mbox{and} \quad
\langle S^* S g,g \rangle + 2 \langle S^*(Sf_\sharp - s + Sh_\sharp), g \rangle \le 0.
$$
Combining these constraints---specifically, multiplying the first by $a$, the second by $b$, and summing---implies that
\begin{align*}
0 & \ge \langle (a R^* R + b S^* S) g,g \rangle + 2 \langle a R^*(Rf_\sharp - r) +
b S^*(S f_\sharp -s) + (a R^* R + b S^* S) h_\sharp, g \rangle\\
& \ge \langle g,g \rangle + 2 \langle h_\sharp, g \rangle,
\end{align*}
where \eqref{Condab1} and \eqref{Condab2} were exploited in the last step.
In other words, one has $0 \ge \|h_\sharp + g \|^2 - \|h_\sharp\|^2$,
i.e.,  $\|h\|^2 \le \|h_\sharp\|^2$,
under the constraints on $h$,
proving that $h_\sharp$ is indeed a maximizer in \eqref{DefHsharp}.
\epf

\subsubsection{The case $R = \Id$}

As mentioned earlier, the case $R = \Id$ corresponds to the choice $\cV = \{0\}$, 
i.e., to a model set $\cK$ being an origin-centered ball in $H$,
and to regularizations being classical Tikhonov regularizations.
The arguments are slightly less involved here than for the case $R \not= \Id$. 
Here is the main result.

\bthm
Let $S$ be an orthogonal projector on $H$ with $\ker(S) \not= \{0\}$
and let $r \in H$, $s \in {\rm ran}(S)$.
The solution $f_{\tau_\sharp}$ to the regularization program \eqref{GenReg2} with parameter
$$
\tau_\sharp = \max \bigg\{ 1-\f{\eta}{\|Sr-s\|},0 \bigg\} 
$$
is the Chebyshev center of the set $\{ f \in H: \| f-r\| \le \eps, \|Sf-s\| \le \eta \}$.
\ethm 

\bpf
Before separating two cases,
we remark that $\|Sr-s\| \le \eps+ \eta$ is implicitly assumed for the above set to be nonempty.
Now,  we first consider the case $\|Sr-s\| > \eta$.
Defining $f_\sharp := f_{\tau_\sharp}$ with $\tau_\sharp = 1-\eta/\|Sr-s\| \in (0,1)$,
our objective is to find $h_\sharp \in H$ and $a,b \ge 0$ for which 
conditions \eqref{CondNorm}, \eqref{Condab1}, and \eqref{Condab2} of Lemma \ref{LemMaximizer} are fulfilled,
so that $h_\sharp$ is a maximizer appearing in Lemma \ref{LemCheCen},
and then to verify that the orthogonality conditions \eqref{CondOrtho} hold,
so that $f_\sharp$ is indeed the required Chebyshev center.
We take any $h_\sharp \in \ker(S)$, with a normalization will be decided later,
and $a = 1$, $b = \tau_\sharp/(1-\tau_\sharp)$. 
In this way,  since $R=\Id$,
condition \eqref{Condab1} is automatic,
and condition~\eqref{Condab2} follows from the characterization~\eqref{CharaGenReg}
written here as $(1-\tau_\sharp)(f_\sharp-r) = -\tau_\sharp (S f_\sharp - s)$.
This characterization also allows us to deduce \eqref{CondOrtho} only from 
$\langle S f_\sharp - s, h_\sharp \rangle = 0$,
which holds because the spaces ${\rm ran}(S)$ and $\ker(S)$ are orthogonal.
The remaining condition \eqref{CondNorm} now reads $\|f_\sharp - r \|^2 + \|h_\sharp\|^2 = \eps^2$ and $\|S f_\sharp - s\|^2 = \eta^2$. 
Recalling from Proposition \ref{PropLinearExpr} that $f_\sharp = (1-\tau_\sharp) f_0 + \tau_\sharp f_1$,
while taking into account that $f_0 = r $ here
and that $f_1 = f_0 + S(f_1-f_0) = r+s-Sr$
thanks to \eqref{f0Nf1},
we have $f_\sharp - r = \tau_\sharp (s-Sr)$
and $S f_\sharp -s = -(1-\tau_\sharp)(s-Sr)$.
Thus, condition \eqref{CondNorm} reads
$$
\tau_\sharp^2 \|s-Sr\|^2 + \|h_\sharp\|^2 = \eps^2
\qquad \mbox{and} \qquad
(1-\tau_\sharp)^2 \|s-Sr\|^2 = \eta^2.
$$ 
The latter is justified by our choice of $\tau_\sharp$,
while the former can simply be achieved by normalizing~$h_\sharp$,
so long as $\eps \ge \tau_\sharp \|s-Sr\|$,
i.e., $\eps \ge \|s-Sr\| - \eta$,
which is our implicit assumption for nonemptiness of the set under consideration. 

Next, we consider the case $\|Sr-s\| \le \eta$.
We note that this implies that $r$ belongs to the set $\{ f \in H: \|f-r\| \le \eps, \|Sf-s\| \le \eta  \}$---we are going to show that $r$ is actually the Chebyshev center of this set.
In other words,  since $r=f_0$, this means that $f_{\tau_\sharp}$ with $\tau_\sharp = 0$ is the Chebyshev center.
To this end, we shall establish that, for any $g \in H$, 
$$
\sup_{\substack{ \|f-r\| \le \eps \\ \|Sf-s\| \le \eta }} \|f-g\|
\ge \sup_{\substack{ \|f-r\| \le \eps \\ \|Sf-s\| \le \eta }} \|f-r\|.
$$
On the one hand,
the right-hand side is obviously bounded above by $\eps$.
On the other hand, selecting $h \in \ker(S)$ with $\|h\| = \eps$,
we define $f_\pm := r \pm h$
to obtain $\|f_\pm - r \| = \|h\| = \eps$
and $\|Sf_\pm -s \| = \|Sr - s\| \le \eta$.
Thus,  the left-hand side is bounded below by
$$
\max_\pm \|f_\pm - g\| \ge \f{1}{2}\|f_+ - g\| + \f{1}{2} \|f_- -g \|
\ge \f{1}{2} \| (f_+ - g) - (f_- - g)\| = \f{1}{2} \|2 h\| = \eps. 
$$
This proves that the left-hand side is larger than or equal to the right-hand side, as required.
\epf

\subsubsection{The case $R \not= \Id$}

We now assume that $R$ is a proper orthogonal projection,
i.e., that $R \not= \Id$,
which corresponds to the case $\cV \not= \{0\}$.
The main result is stated below.
To apply it in practice,
the optimal parameter $\tau$ needs to be computed by solving an equation involving the smallest eigenvalue of self-adjoint operator depending on $\tau$.
This can be done using an all purpose routine.
We could also devise our own bisection method,
%(exploiting the uniqueness of the Chebyshev center, and hence of the optimal $\tau$),
Newton method (since the derivative $d\la_{\min}/d\tau$ is accessible, see the appendix),
or secant method.

\bthm
\label{ThmLocalTF}
Let $R \not= \Id, S \not= \Id$ be two orthogonal projectors on $H$ such that $\ker(R) \cap \ker(S) = \{0\}$
and let $r \in {\rm ran}(R)$, $s \in {\rm ran}(S)$.
Consider $\tau_\sharp$ to be a (often unique) $\tau$ between $1/2$ and $\eps/(\eps+\eta)$ such that
\label{PageEigEq}
\be
\label{EqForLocal}
\la_{\min}((1-\tau) R + \tau S) - \f{(1-\tau)^2 \eps^2 - \tau^2 \eta^2}{(1-\tau) \eps^2 - \tau \eta^2 + (1-\tau)\tau(1-2\tau) \delta^2}
=0,
\ee
where $\delta$ is precomputed as $\delta = \min\{ \|Rf-r\|: Sf=s\} = \min\{ \|Sf-s\|: Rf=r \}$.
Then the solution $f_{\tau_\sharp}$ of the regularization program \eqref{GenReg2} with parameter $\tau_\sharp$
is the Chebyshev center of the set $\{ f \in H: \|Rf-r\|\le \eps, \|Sf-s\| \le \eta \}$.
\ethm

\brk
It there is no observation error, 
i.e., if $\eta = 0$,
then the parameter solving equation \eqref{EqForLocal} is $\tau_{\sharp} = 1$.
In case $R = P_{\cV^\perp}$, $r=0$,  $S = \La^* \La$, and $s= \La^* y$,
this means that the Chebyshev center is $f_1 = \argmin \|P_{\cV^\perp} f\|$ s.to $\La f = y$
and we thus retrieve the result of \citep*{binev2017data}.
\erk

%\brk
%We can probably deal with the recovery of a quantity of interest $Q(f)$ instead of just $f$...
%\erk

The proof of Theorem \ref{ThmLocalTF} requires an additional result that gives information about 
the norms of the projections $Rh$ and $Sh$  
when $h$ is  an eigenvector of the positive semidefinite operator  $(1-\tau)R + \tau S$.
This result will be applied for the eigenvector associated with the smallest eigenvalue.

\blem
\label{LemEigenStuff}
Let $R,S$ be two orthogonal projectors on $H$.
For $\tau \in (0,1)$,
let $h \in H$ be an eigenvector of $(1-\tau) R + \tau S$ corresponding to an eigenvalue $\la \not= 1/2$.
Then
\be
\label{RhSh}
 \|Rh\|^2 = \f{(\tau-\la)\la}{(1-\tau)(1-2\la)} \|h\|^2
\qquad \mbox{and} \qquad
\|Sh\|^2 =  \f{(1-\tau-\la)\la}{\tau (1-2\la)} \|h\|^2.
\ee
\elem

\bpf
We notice, on the one hand, that
\begin{align}
\label{EstRhSh1}
(1-\tau) \|Rh\|^2 + \tau \|Sh\|^2 
& = (1-\tau) \langle Rh,h \rangle + \tau \langle Sh,h \rangle = \langle ( (1-\tau) R + \tau  S)h,h \rangle\\
\nonumber
& = \la \|h\|^2,
\end{align}
and, on the other hand,  that
\begin{align*}
(1-\tau)^2 \|Rh\|^2 - \tau^2 \|Sh\|^2 
& = \langle (1-\tau) Rh + \tau Sh,  (1-\tau) Rh - \tau S h\rangle
= \langle \la h,  (1-\tau) Rh - \tau S h\rangle\\
& = \la (1-\tau) \|Rh\|^2 - \la \tau \|Sh\|^2.
\end{align*}
Rearranging the latter yields
\be
\label{EstRhSh2}
(1-\tau)(1-\tau-\la) \|Rh\|^2 -  \tau(\tau-\la)\|Sh\|^2 = 0.
\ee
Together,  the equaltions \eqref{EstRhSh1} and \eqref{EstRhSh2} form a two-by-two linear system in the unknowns $\|Rh\|^2$ and $\|Sh\|^2$ with determinant  $-(1-\tau)\tau(1-2\la) \not= 0$.
Its solutions are easily verified to be the ones given in \eqref{RhSh}.
\epf

\brk
Because $\|Rh\|^2$, $\|Sh\|^2$, and $\|h\|^2$ are all nonnegative,
Lemma \ref{LemEigenStuff} implicitly guarantees that $\tau - \la$ and $1-\tau-\la$ have the same sign as $1-2 \la \not=0$.
These quantities are nonnegative when $R \not= \Id$, $S \not= \Id$,  and $\la$ is the smallest eigenvalue---the case of application of the lemma.
Indeed,  taking $f \in \ker(R)$ with $\|f\|=1$ (which is possible because $R \not= \Id$),
one has
$$
\la_{\min}:=\la_{\min}((1-\tau)R+\tau S)
\le \|(1-\tau) Rf + \tau Sf \| = \tau \|Sf\| \le \tau,
$$
i.e., $\tau - \la_{\min} \ge 0$.
The inequality  $\la_{\min} \le 1-\tau$, i.e., $1-\tau-\la_{\min} \ge 0$,  is obtained in a similar fashion.
These inequalities sum up to give $1-2 \la_{\min} \ge 0$.
The latter is in fact (strictly) positive when $\tau \not= 1/2$,
since either $\tau$ or $1-\tau$ is smaller than $1/2$, so that $\la_{\min} < 1/2$.
\erk

%\brk
%In the common situation $\cV = \{0\}$,
%we take $R = P_{\cV^\perp} = \Id$,
%so $\la_{\min}((1-\tau) R + \tau S) = 1-\tau + \la_{\min}(S) = 1-\tau$
%and a corresponding eigenvector satisfies $Sh=0$.
%Note that \eqref{RhSh} are still valid in this case with $\|Rh\|^2 = \|h\|^2$ and $\|Sh\|^2 = 0$.
%\erk

With the above result at hand,
we are ready to fully justify the main result of this subsection.

\bpf[Proof of Theorem \ref{ThmLocalTF}]
Let us temporarily take for granted the existence of a solution $\tau_\sharp$ to \eqref{EqForLocal}.
Defining $f_\sharp := f_{\tau_\sharp}$,
our objective is again to find $h_\sharp \in H$ and $a,b \ge 0$ for which 
conditions \eqref{CondNorm}, \eqref{Condab1}, and \eqref{Condab2} of Lemma \ref{LemMaximizer} are fulfilled,
so that $h_\sharp$ is a maximizer appearing in Lemma \ref{LemCheCen},
and then to verify that the orthogonality conditions \eqref{CondOrtho} hold,
so that $f_\sharp$ is indeed the required Chebyshev center.
Writing $\la_\sharp := \la_{\min}((1-\tau_\sharp) R + \tau_\sharp S)$,
we choose $h_\sharp$ to be a (so far unnormalized) eigenvector of $(1-\tau_\sharp) R + \tau_\sharp S$ corresponding to the eigenvalue $\la_\sharp$.
Setting $a:= (1-\tau_\sharp) / \la_\sharp$ and $b:= \tau_\sharp/\la_\sharp$,
conditions \eqref{Condab1}  is swiftly verified, since
$RR^* = R$, $SS^*=S$, and
$$
a R + b S = \f{(1-\tau_\sharp)  R + \tau_\sharp  S}{\la_{\min}((1-\tau_\sharp) R + \tau_\sharp S)}
\succeq \Id.
$$
Then,  the characterization $(1-\tau_\sharp) R(f_\sharp - r) = -\tau_\sharp S(f_\sharp - s)$ of the regularization solution $f_\sharp$,
see~\eqref{CharaGenReg},
allows us to validate condition \eqref{Condab2} via
\begin{align*}
a R(f_\sharp - r) + b S(f_\sharp - s) + (a R + b S ) h_\sharp 
& = \f{1}{\la_\sharp}
\Big( (1-\tau_\sharp) R(f_\sharp - r) \hspace{-.5mm}+\hspace{-.5mm} \tau_\sharp S(f_\sharp - s)
\hspace{-.5mm}+\hspace{-.5mm} ((1-\tau_\sharp)  R + \tau_\sharp S) h_\sharp  \Big)\\
& =\f{1}{\la_\sharp} \left( 0 + \la_\sharp h_\sharp \right) = h_\sharp.
\end{align*}
The orthogonality conditions \eqref{CondOrtho} are also swiftly verified: 
the second one follows from the first one using \eqref{CharaGenReg};
the first one holds because, 
while $h_\sharp$ is an eigenvector of $(1-\tau_\sharp) R + \tau_\sharp S$ corresponding to its smallest eigenvalue,
$R (f_{\sharp} - r)  = -\tau_\sharp/(1-\tau_\sharp) S(f_\sharp - s)$ is an eigenvector corresponding to the largest eigenvalue (i.e., to one),
since it is invariant when applying both $R$ and $S$. 
Thus, it remains to verify that the two conditions of \eqref{CondNorm} are fulfilled.
In view of the orthogonality conditions \eqref{CondOrtho},
they read
\be
\label{Identities4Later}
\|R f_\sharp - r \|^2 + \|R h_\sharp\|^2 = \eps^2
\qquad \mbox{and} \qquad
\|S f_\sharp - s \|^2 + \|S h_\sharp\|^2  = \eta^2.
\ee
Now, invoking Proposition \ref{PropLinearExpr},
as well as Lemma \ref{LemEigenStuff}, the two conditions of \eqref{CondNorm}  become
\begin{align}
\label{NeededCond1}
\tau_\sharp^2 \delta^2 + \f{(\tau_\sharp-\la_\sharp)\la_\sharp}{(1-\tau_\sharp) (1-2 \la_\sharp)} \|h_\sharp\|^2 & = \eps^2\\
\label{NeededCond2}
(1-\tau_\sharp)^2 \delta^2 + \f{(1-\tau_\sharp-\la_\sharp)\la_\sharp}{\tau_\sharp (1-2 \la_\sharp)} \|h_\sharp\|^2 & = \eta^2.
\end{align}
After some simplification work,
starting by forming the combinations $(1-\tau_\sharp)\times$\eqref{NeededCond1}$-\tau_\sharp^2\times$\eqref{NeededCond2}
and $(1-\tau_\sharp-\la_\sharp)(1-\tau_\sharp)\times$\eqref{NeededCond1}$-(\tau_\sharp-\la_\sharp)(\tau_\sharp)\times$\eqref{NeededCond2},
these two conditions are seen to be equivalent to
\begin{align}
\label{ThisIsLWCE}
\|h_\sharp\|^2 & = 
\f{1-2\la_\sharp}{(2\tau_\sharp - 1)\la_\sharp^2}\big( (1-\tau_\sharp)^2 \eps^2
- \tau_\sharp^2 \eta^2  \big),\\
\label{ThisIsLaMin}
\la_{\sharp} & = 
\f{ (1-\tau_\sharp)^2 \eps^2 - \tau_\sharp^2 \eta^2 }{ (1-\tau_\sharp) \eps^2 
- \tau_\sharp \eta^2 + (1-\tau_\sharp) \tau_\sharp (1- 2\tau_\sharp) \delta^2}.
\end{align}
These two conditions can be fulfilled:
the latter is the condition that defined $\tau_\sharp$,
i.e., \eqref{EqForLocal},
while the former is simply guaranteed by properly normalizing the eigenvector $h_\sharp$.

Before establishing the existence $\tau_\sharp$,
we point out that its uniqueness holds when $f_0 \not= f_1$,
i.e., when there is no $f \in H$ such that $Rf=r$ and $Sf=s$---such an $f$ would solve the regularization program for any $\tau \in [0,1]$.
Indeed,  if $\tau \not= \tau'$ were two solutions to \eqref{EqForLocal},
then the previous argument would imply that $f_\tau$ and $f_{\tau'}$ are both Chebyshev centers,
which could only happen if they were equal, 
i.e., if $f_0=f_1$ by \eqref{ExprGenRegAff}.
Now,  for the existence of $\tau_\sharp$,
it will be justified by the fact that the function  
$$
\theta: \tau \mapsto \la_{\min}((1-\tau) R + \tau S) - \f{(1-\tau)^2 \eps^2 - \tau^2 \eta^2}{(1-\tau) \eps^2 - \tau \eta^2 + (1-\tau)\tau(1-2\tau) \delta^2}
$$
is continuous between $1/2$
and $\eps / (\eps + \eta)$ and takes values of different signs there.
To see the difference in sign,  notice that $\la_{\min}((1-\tau) R + \tau S) \in [0,1/2]$ by the remark after Lemma \ref{LemEigenStuff}---this is where the assumption $R \not= \Id$ is critical---so that
$$
\theta\bigg( \f{1}{2} \bigg) \le \f{1}{2} - \f{1}{2} \le0
\qquad \mbox{and} \qquad 
\theta \bigg( \f{\eps}{\eps+\eta} \bigg) \ge 0  - 0 \ge 0.
$$
To see the continuity,
we need the continuity of the smallest eigenvalue as a function of $\tau$
and the nonvanishing of the denominator $(1-\tau) \eps^2 - \tau \eta^2 + (1-\tau)\tau(1-2\tau) \delta^2$ between $1/2$ and $\eps/(\eps+\eta)$.
The former is a consequence of Weyl's inequality, yielding 
$$
| \la_{\min}((1-\tau) R + \tau S) - \la_{\min}((1-\tau') R + \tau' S) |
 \le \| ((1-\tau) R + \tau S) - ((1-\tau') R + \tau' S) \| 
 = |\tau - \tau'| \, \|R-S\|.
$$
The latter is less immediate.
We start by using \eqref{NormRFtau-r} and recalling the very definition of $f_\tau$ to write
$$
(1-\tau) \tau \delta^2
= (1-\tau) \|Rf_\tau - r\|^2 + \tau \|Sf_\tau -r \|^2
\le (1-\tau) \eps^2 + \tau \eta^2.
$$
Therefore,  if the denominator vanished for some $\tau \in (0,1) \setm \{1/2\}$,
we would have
\begin{align*}
0 & = \f{(1-\tau) \eps^2 - \tau \eta^2}{1-2\tau} + (1-\tau) \tau \delta^2
\le \f{(1-\tau) \eps^2 - \tau \eta^2}{1-2\tau} + (1-\tau) \eps^2 + \tau \eta^2
= \f{2(1-\tau)^2 \eps^2 - 2 \tau^2 \eta^2}{1-2\tau}\\
& = \f{ \big( (1-\tau)\eps + \tau \eta \big)  \big((1-\tau)\eps - \tau \eta \big)}{1/2-\tau}
=  \big( (1-\tau)\eps + \tau \eta \big) \big( \eps + \eta \big) \f{\eps/(\eps+\eta)-\tau}{1/2-\tau}.
\end{align*}
This would force $\eps/(\eps+\eta) - \tau$ and $1/2 - \tau$ to have the same sign,
contrary to the assumption that $\tau$ runs between $1/2$ and $\eps/(\eps+\eta)$.
Thus, the nonvanishing of the denominator is explained,
concluding the proof.
\epf

\brk
The above arguments contain the value of the minimal local worst-case error,
i.e., of the Chebyshev radius of the set $\cC = \{ f \in H: \|Rf-r\| \le \eps,  \|Sf-s\| \le \eta \}$.
Indeed, we recall from the proof of Lemma \ref{LemCheCen} that this radius equals $\|h_\sharp\|$, whose value was derived in \eqref{ThisIsLWCE}.
This expression can be simplified with the help of \eqref{ThisIsLaMin} by noticing that
$$
\frac{1-2\la_\sharp}{\la_\sharp} = (2\tau_\sharp - 1) \f{(1-\tau_\sharp) \eps^2 + \tau_\sharp \eta^2 -(1-\tau_\sharp) \tau_\sharp \delta^2}{(1-\tau_\sharp)^2 \eps^2 - \tau_\sharp^2 \eta^2}.
$$
As a consequence, we deduce that the Chebyshev radius satisfies
$$
{\rm radius}(\cC)^2
= \f{1-\tau_\sharp}{\la_\sharp} \eps^2 + \f{\tau_\sharp}{\la_\sharp} \eta^2
- \f{(1-\tau_\sharp)\tau_\sharp}{\la_\sharp} \delta^2,
\qquad \quad
\la_\sharp:= \la_{\min}((1-\tau_\sharp)R+\tau_\sharp S).
$$
\erk

\section{Global Optimality}
\label{SecGlob}

Our goal in this section is to uncover some favorable globally optimal recovery maps---favorable in the sense that they are linear maps.
We start by considering the situation of an arbitrary observation map $\La$
before moving to the particular case where it satisfies $\La \La^* = \Id_{\bR^m}$.

\subsection{Arbitrary observations}

In this subsection,
we first recall
a standard lower bound for the global worst-case error. 
This lower bound,  already exploited e.g. in \citep{micchelli1993optimal},
shall be expressed
%as well as the (squared) global worst-case error of a linear recovery map,
as the minimal value of a certain semidefinite program.
This expression will allow us to demonstrate that the lower bound is achieved by the regularization map
$$
\Delta_\tau: y \in \bR^m \mapsto
\underset{f \in H}{\argmin} \; (1-\tau) \| P_{\cV^\perp} f\|^2 + \tau \|\La f - y\|^2 
$$
for some parameter $\tau \in (0,1)$ to be explicitly determined.
Here is a precise formulation of the result.

\bthm
\label{ThmGlobalGen}
Given the approximability set $\cK = \{f \in H: {\rm dist}(f,\cV) \le \eps \}$
and the uncertainty set $\cE = \{e \in \bR^m: \|e\|\le \eta \}$,
define $\tau_\flat := d_\flat/(c_\flat+d_\flat)$ where
$c_\flat,d_\flat \ge 0$ are solutions to 
$$
\minimize{c,d \ge 0} \; c \eps^2 + d \eta^2
\qquad \mbox{s.to} \quad
c P_{\cV^\perp} + d \La^* \La  \succeq \Id.
$$
Then the regularization map $\Delta_{\tau_\flat}$ is a globally optimal recovery map over $\cK$ and $\cE$, i.e.,
\be
{\rm gwce}(\Delta_{\tau_\flat}) = \inf_{\Delta: \bR^m \to H}
{\rm gwce}(\Delta).
\ee
\ethm

The proof relies on three lemmas given below,  the first of which introducing the said lower bound.

\blem
\label{LemIsolated1}
For any recovery map $\Delta: \bR^m \to H$, one has ${\rm gwce}(\Delta) \ge {\rm lb}$, where
$$
%\label{TheLB}
{\rm lb} := \sup_{\substack{\|P_{\cV^\perp}h\| \le \eps\\ \|\La h\|\le \eta}} \|h\|.
$$
\elem

\bpf
As a reminder, the global worst-case error of $\Delta$ is defined by  
$$
{\rm gwce}(\Delta) = \sup_{\substack{ \|P_{\cV^\perp} f\| \le \eps \\ \|e\| \le \eta }}
\| f - \Delta(\La f+e) \|.
$$
For any $h \in H$ such that $\|P_{\cV^\perp}h\| \le \eps$ and $\| \La h \| \le \eta $,
since $f_\pm = \pm h$ satisfies $\|P_{\cV^\perp} f_\pm \| \le \eps$
and $e_\pm = \mp \La h $ satisfies $\|e_\pm\| \le \eta$, we have
\begin{align*}
{\rm gwce}(\Delta) & \ge \max_\pm \| f_\pm - \Delta( \La f_\pm + e_\pm )\|
= \max_\pm \| f_\pm - \Delta(0) \| 
 \ge \f{1}{2} \|f_+ - \Delta(0) \| + \f{1}{2} \|f_- - \Delta(0) \|\\
& \ge \f{1}{2} \| (f_+-\Delta(0)) - (f_- - \Delta(0))  \|
= \f{1}{2} \|2 h\| = \|h\|.
\end{align*}
Taking the supremum over $h$ leads to the required inequality ${\rm gwce}(\Delta) \ge {\rm lb}$.
\epf

The second lemma expresses the square of the lower bound as the minimal value of a semidefinite program.
In passing, 
the square of the global worst-case error of a linear recovery map
is also related to the minimal value of a semidefinite program.

\blem
\label{LemIsolated2}
One has
\be
\label{LB4GWCE}
{\rm lb}^2  = \min_{c,d \ge 0} \; c \eps^2 + d \eta^2
\qquad \mbox{s.to} \quad
c P_{\cV^\perp} + d \La^* \La  \succeq \Id.
\ee
Moreover, if a recovery map $\Delta: \bR^m \to H$ is linear, one also has
\be
\label{GWCETau}
{\rm gwce}(\Delta)^2  \le 
\min_{c,d \ge 0} \; c \eps^2 + d \eta^2
\qquad \mbox{s.to} \quad
\bbmx
c P_{\cV^\perp} & | & 0\\
\hline
0 & | & d \, \Id_{\bR^m}
\ebmx
\succeq \bbmx \Id - \La^* \Delta^*\\ \hline  \Delta^* \ebmx
\bbmx \Id - \Delta \La \; | \; \Delta \ebmx. 
\ee
\elem

\bpf
The first semidefinite characterization is based on the version of the S-procedure stated in Theorem \ref{ThmPolyak}.
Precisely, we write the square of the lower bound as
\begin{align*}
{\rm lb}^2 
& =\;\,  \sup_h \, \|h\|^2 & & \mbox{s.to }
\|P_{\cV^\perp} h\|^2 \le \eps^2 \mbox{ and } \|\La h \|^2 \le \eta^2\\
 & = \; \,\inf_{\gamma} \, \gamma
& & \mbox{s.to } \|h\|^2 \le \gamma \mbox{ whenever } \|P_{\cV^\perp} h \|^2 \le \eps^2 
\mbox{ and } \|\La h \|^2 \le \eta^2\\
& = \; \,\inf_{\gamma} \, \gamma& & \mbox{s.to } \exists \, c,d \ge 0: \;
 \|h\|^2 - \gamma \le c(\|P_{\cV^\perp} h \|^2 - \eps^2)  + d(\|\La h \|^2 - \eta^2)
 \mbox{ for all }h \in H\\
 & = \inf_{\substack{\gamma \\ c,d \ge 0}} \, \gamma & &\mbox{s.to } 
c  \langle P_{\cV^\perp} h, h \rangle + d \langle \La^* \La h, h \rangle - \langle h,h \rangle + \gamma - c \eps^2 - d \eta^2 \ge 0   \mbox{ for all }h \in H.
\end{align*}
The validity of Theorem \ref{ThmPolyak} in ensured by the facts that $\|P_{\cV^\perp}\wt{h}\|^2 - \eps^2 < 0$
and $\|\La \wt{h}\|^2 - \eta^2 < 0$ for $\wt{h} = 0$
and that $P_{\cV^\perp} + \La^* \La \succ 0$.
Note that the resulting constraint decouples as
$\langle c P_{\cV^\perp}h  +  d \La^* \La h- h,h \rangle \ge 0$ for all $h \in H$,
i.e.,  $cP_{\cV^\perp} + d \La^* \La  - \Id \succeq 0$,
and $\gamma - c \eps^2 - d \eta^2 \ge 0$.
Taking the minimal value of $\gamma$ under the latter constraint, 
namely $c \eps^2 + d \eta^2$,
leads to the expression of ${\rm lb}^2$ given in \eqref{LB4GWCE}.

As for \eqref{GWCETau},
we start by remarking that the linearity of the recovery map $\Delta$ allows us to write
\begin{align*}
{\rm gwce}(\Delta)^2 
& = \sup_{f,e} \, \|f - \Delta \La f - \Delta e\|^2 \quad \mbox{s.to }
\|P_{\cV^\perp} f\|^2 \le \eps^2 \mbox{ and } \|e \|^2 \le \eta^2\\
& = \inf_\gamma \, \gamma \quad \mbox{s.to }
\|f - \Delta \La f - \Delta e\|^2 \le \gamma  \mbox{ whenever } 
\|P_{\cV^\perp} f\|^2 \le \eps^2 \mbox{ and } \|e \|^2 \le \eta^2.
\end{align*}
The latter constraint can be expressed in terms of the combined variable $v=(f,-e) \in H \times \bR^m$ as
\be
\label{CstForGWCE}
\Big\| \bbmx  \Id - \Delta \La \;  |  \; \Delta \ebmx v \Big\|^2 \le \gamma  \mbox{ whenever } 
\Big\| \bbmx P_{\cV^\perp} \; | \; 0 \ebmx v \Big\|^2 \le \eps^2 \mbox{ and }
\Big\| \bbmx 0 \; | \; \Id_{\bR^m} \ebmx v \Big\|^2 \le \eta^2.
\ee
Although the proviso of Theorem \ref{ThmPolyak} is not fulfiled here,
the constraint \eqref{CstForGWCE} is still a consequence of (but is not equivalent to) the existence of $c,d \ge 0$ such that
$$
\; \Big\| \bbmx  \Id - \Delta \La \; | \; \Delta \ebmx v \Big\|^2 - \gamma
\le c \Big(  \Big\| \bbmx P_{\cV^\perp} \; | \; 0 \ebmx v \Big\|^2 - \eps^2 \Big)   + d \Big( \Big\| \bbmx 0 \; | \; \Id_{\bR^m} \ebmx v \Big\|^2 - \eta^2 \Big)
\quad
 \mbox{for all } v \in H \times \bR^m.
$$
The latter can also be written as the existence of $c,d \ge 0$ such that,
for all  $v \in H \times \bR^m$,
\begin{multline*}
\Big\langle  \Big( c \bbmx P_{\cV^\perp} \; | \; 0 \ebmx^* \bbmx P_{\cV^\perp} \; | \; 0 \ebmx + d \bbmx 0 \; | \; \Id_{\bR^m} \ebmx^* \bbmx 0 \; | \; \Id_{\bR^m} \ebmx - \bbmx  \Id - \Delta \La \; | \; \Delta \ebmx^* \bbmx  \Id - \Delta \La \; | \; \Delta \ebmx \Big) v  ,v \Big\rangle \\
+ \gamma - c \eps^2 - d \eta^2 \ge 0.
\end{multline*}
Therefore, we obtain the inequality (instead of the equality)
\begin{align*}
{\rm gwce}(\Delta)^2
& \le \inf_{\substack{\gamma\\c,d \ge 0}} \, \gamma
& \mbox{s.to} \quad & 
c\bbmx
P_{\cV^\perp} & | & 0\\
\hline
0 & | & 0
\ebmx
+
d \bbmx
0 & | & 0\\
\hline
0 & | &  \Id_{\bR^m}
\ebmx
- 
\bbmx \Id - \La^* \Delta^*\\ \hline  \Delta^* \ebmx
\bbmx \Id - \Delta \La \; | \; \Delta \ebmx \succeq 0\\
& & \mbox{and} \quad & \gamma - c \eps^2 - d \eta^2 \ge 0.
\end{align*}
The variable $\gamma$ can be eliminated from this optimization program
by assigning it the value $c \eps^2 + d \eta^2$,
thus arriving at the semidefinite program announced in \eqref{GWCETau}.
\epf

The third and final lemma relates the constraints of \eqref{LB4GWCE}  and \eqref{GWCETau}:
while the constraint of \eqref{GWCETau} with any regularization map $\Delta_\tau$ implies the constraint of \eqref{LB4GWCE}, see the appendix,
we need the partial converse that the constraint of \eqref{LB4GWCE} 
implies the constraint of \eqref{GWCETau} for a specific regularization map~$\Delta_\tau$.
\label{PageImplies1}

\blem
\label{LemIsolated3}
If $c P_{\cV^\perp} + d \La^* \La \succeq \Id$,
then setting $\tau = d/(c+d)$ yields
$$
\bbmx
c P_{\cV^\perp} & | & 0\\
\hline
0 & | & d \,  \Id_{\bR^m}
\ebmx 
\succeq
 \bbmx \Id - \La^* \Delta_\tau^*\\ \hline  \Delta_\tau^* \ebmx
\bbmx \Id - \Delta_\tau \La \; | \; \Delta_\tau \ebmx. 
$$
\elem

\bpf
We recall from Proposition \ref{PropGenExpr} adapted to the current situation that, for any $\tau \in (0,1)$,
$$
\Delta_\tau = \big( (1-\tau) P_{\cV^\perp} + \tau \La^* \La \big)^{-1} (\tau \La^*),
\quad \mbox{hence} \quad
\Id - \Delta_\tau \La = \big( (1-\tau) P_{\cV^\perp} + \tau \La^* \La \big)^{-1} ((1-\tau) P_{\cV^\perp}).
$$ 
We now  notice that the hypothesis  $c P_{\cV^\perp} + d \La^* \La \succeq \Id$ is equivalent to $\la_{\min} (c P_{\cV^\perp} + d \La^* \La) \ge 1$.
With our particular choice of $\tau$, this reads
$\la_{\min} ( (1-\tau) P_{\cV^\perp} + \tau \La^* \La ) \ge 1/(c+d)$.
It  follows that 
$$
\la_{\max} \big( ((1-\tau)P_{\cV^\perp} + \tau \La^* \La )^{-1}  \big) 
= \f{1}{\la_{\min}((1-\tau)P_{\cV^\perp} + \tau \La^* \La)} 
\le c+d.
$$
The inverse appearing above can be written as  
$$
\big( (1-\tau) P_{\cV^\perp} + \tau \La^* \La \big)^{-1} \bbmx  \sqrt{1-\tau} P_{\cV^\perp} \; | \; \sqrt{\tau} \La^* \ebmx
\bbmx  \sqrt{1-\tau} P_{\cV^\perp} \\ \hline \sqrt{\tau} \La \ebmx \big( (1-\tau) P_{\cV^\perp} + \tau \La^* \La \big)^{-1} 
$$
and since $AB$ and $BA$ always have the same nonzero eigenvalues, we derive that
$$
\la_{\max} \left(
\bbmx  \sqrt{1-\tau} P_{\cV^\perp} \\ \hline \sqrt{\tau} \La \ebmx
\big( (1-\tau) P_{\cV^\perp} + \tau \La^* \La  \big)^{-2}
\bbmx  \sqrt{1-\tau} P_{\cV^\perp} \; | \; \sqrt{\tau} \La^* \ebmx
\right) 
\le c+d.
$$
Writing the latter as
$$
\bbmx  \sqrt{1-\tau} P_{\cV^\perp} \\ \hline \sqrt{\tau} \La \ebmx
\big( (1-\tau) P_{\cV^\perp} + \tau \La^* \La  \big)^{-2}
\bbmx  \sqrt{1-\tau} P_{\cV^\perp} \; | \; \sqrt{\tau} \La^* \ebmx 
\preceq (c+d) \Id
$$
and multiplying on both sides by 
$\bbmx \sqrt{1-\tau} P_{\cV^\perp} & | & 0\\ \hline 0 & | & \sqrt{\tau} \, \Id_{\bR^m} \ebmx$
yields
$$
\bbmx  (1-\tau) P_{\cV^\perp} \\ \hline \tau \La \ebmx
 \big( (1-\tau) P_{\cV^\perp} +  \tau \La^* \La  \big)^{-2}
\bbmx  (1-\tau) P_{\cV^\perp} \; | \; \tau \La^* \ebmx 
\preceq (c+d) 
\bbmx (1-\tau) P_{\cV^\perp} & | & 0\\ \hline 0 & | & \tau \, \Id_{\bR^m} \ebmx.
$$
Taking the expressions of $\Delta_\tau$ and $\Id - \Delta_\tau \La $ into account, we conclude that
$$
 \bbmx \Id - \La^* \Delta_\tau^*\\ \hline  \Delta_\tau^* \ebmx
\bbmx \Id - \Delta_\tau \La \; | \; \Delta_\tau \ebmx
\preceq 
\bbmx
c P_{\cV^\perp} & | & 0\\
\hline
0 & | & d \,  \Id_{\bR^m}
\ebmx,
$$
as announced.
\epf

With the above three lemmas at hand,
the main result of this subsection follows easily.

\bpf[Proof of Theorem \ref{ThmGlobalGen}]
Since Lemma \ref{LemIsolated1} guarantees that 
$ \inf \{ {\rm gwce}(\Delta), \Delta: \bR^m \to H \} \ge {\rm lb}$,
we only need to show that  ${\rm gwce}(\Delta_{\tau_\flat}) \le {\rm lb}$.
By the first part of Lemma \ref{LemIsolated2},
we have ${\rm lb}^2 = c_\flat \eps^2 + d_\flat \eta^2$
with $c_\flat$ and $d_\flat$ satisfying $c_\flat P_{\cV^\perp} + d_\flat \La \La^* \succeq \Id$.
By Lemma \ref{LemIsolated3}, the latter implies that
$$
\bbmx
c_\flat P_{\cV^\perp} & | & 0\\
\hline
0 & | & d_\flat \Id_{\bR^m}
\ebmx 
\succeq
 \bbmx \Id - \La^* \Delta_{\tau_\flat}^*\\ \hline  \Delta_{\tau_\flat}^* \ebmx
\bbmx \Id - \Delta_{\tau_\flat} \La \; | \; \Delta_{\tau_\flat} \ebmx. 
$$
By the second part of Lemma \ref{LemIsolated2},
it follows that ${\rm gwce}(\Delta_{\tau_\flat})^2 \le c_\flat \eps^2 + d_\flat \eta^2 = {\rm lb}^2$, which is the required inequality.\epf

\brk
When $\cV = \{0\}$, so that $P_{\cV^\perp} = \Id$,
we obtain $c_\flat = 1$ and  $d_\flat =0$,
resulting in a minimal global worst-case error equal to $\eps$
and achieved for the regularization map $\Delta_0=0$.
This result can be seen directly from
${\rm gwce}(\Delta) \ge \sup\{ \|h\|: \|h\|\le \eps, \|\La h \| \le \eta \} = \eps$
for any $\Delta: \bR^m \to H$, 
while ${\rm gwce}(\Delta_0) = \sup\{ \|f\|: \|f\|\le \eps\} = \eps$.
\erk

%\bpf
%Let $c_*,d_* \ge 0$ be the optimizers of the lower-bound SDP \eqref{LB}.
%We prove below that
%\be 
%\tau_* = \f{d_*}{c_*+d_*}
%\ee
%is an optimal regularization parameter.
%Noticing that $1-\tau_* = c_*/(c_*+d_*)$,
%we first claim that $c_*+d_* =1/( \la_{\min}((1-\tau_*)P_{\cV^\perp} + \tau_* \La^* \La))$.
%Indeed, the feasibility of $c_*,d_*$ reads $(1-\tau_*)P_{\cV^\perp} + \tau_* \La \La^* \succeq 1/(c_*+d_*) \Id$,  so that $\la_{\min}((1-\tau_*)P_{\cV^\perp} + \tau_* \La^* \La) \ge 1/(c_*+d_*)$.
%Now suppose that $\la_{\min}((1-\tau_*)P_{\cV^\perp} + \tau_* \La^* \La) > 1/(c_*+d_*)$.
%Notice that $c = (1-\tau_*)/\la_{\min}, d =\tau_*/\la_{\min}$ are feasible for \eqref{LB},
%so a contradiction arises from
%\be
%c_* \eps^2 + d_* \eta^2 \le \f{1-\tau_*}{\la_{\min}} \eps^2 + \f{\tau_*}{\la_{\min}} \eta^2
%< (1-\tau_*)(c_*+d_*) \eps^2 + \tau_* (c_*+d_*)  \eta^2
%= c_* \eps^2 + d_* \eta^2.
%\ee
%
%Next, we claim that $c_*,d_*$ satisfy the constraint of \eqref{CstTau} for the choice $\tau = \tau_*$.
%This implies that ${\rm gwce}^2(\Delta^{(\tau_*)}) \le c_* \eps^2 + d_* \eta^2$, which is the lower bound, so the announced result is established.
%To justify our claim, we start from
%\epf
%

\subsection{Orthonormal observations}

In this subsection, we demonstrate that
the use of orthonormal observations
guarantees,  rather unexpectedly,  
that regularization provides optimal recovery maps
even without a careful parameter selection.
The main result reads as follows.

\bthm
\label{ThmIndepTau}
Given the approximability set $\cK = \{f \in H: {\rm dist}(f,\cV) \le \eps \}$
and the uncertainty set $\cE = \{e \in \bR^m: \|e\|\le \eta \}$,
if $\La \La^* = \Id_{\bR^m}$,
then all the regularization maps $\Delta_\tau$ are optimal recovery maps, i.e.,
for all $\tau \in [0,1]$,
\be
{\rm gwce}(\Delta_{\tau}) = 
\inf_{\Delta: \bR^m \to H}
{\rm gwce}(\Delta).
\ee
\ethm

\label{PageImplies2}
The proof strategy consists in establishing that the constraints in 
\eqref{LB4GWCE} and in \eqref{GWCETau} with $\Delta = \Delta_\tau$ 
are in fact equivalent for any $\tau \in [0,1]$.
This yields the inequality ${\rm gwce}(\Delta_{\tau}) \le {\rm lb}$,
which proves the required result, given that ${\rm lb}$ was introduced as a lower bound on ${\rm gwce}(\Delta)$ for every $\Delta$.
While the constraint in \eqref{GWCETau} implies the constraint in \eqref{LB4GWCE} for any observation map $\La$  (see the appendix),
the reverse implication relies on the fact that $\La \La^* = \Id_{\bR^m}$,
e.g. via the identity $\Delta_\tau = (1-\tau)\Delta_0 + \tau \Delta_1$ derived in Proposition \ref{PropLinearExpr}.
The following realization is also a crucial point of our argument.

\blem
\label{LemEigcd}
Assume that $\La \La^* = \Id_{\bR^m}$.
For $c,d \ge 0$, 
let $h$ be an eigenvector of $c P_{\cV^\perp} + d \La^* \La$
associated with an eigenvalue $\la$.
For any $\tau \in [0,1]$, one has\vspace{-5mm}
\begin{itemize}
\item if $\la \not= c+d$, then \vspace{-3mm}
$$
(\Id-\La^* \Delta_\tau^*) h = \df{c}{\la} P_{\cV^\perp} h
\qquad \mbox{and} \qquad 
\La^* \Delta_\tau^* h = \df{d}{\la} \La^* \La h;
$$
\item if $\la = c+d$, then  \vspace{-3mm}
$$
(\Id-\La^* \Delta_\tau^*) h =(1-\tau) h
\qquad \mbox{and} \qquad 
\La^* \Delta_\tau^* h = \tau h.
$$
\end{itemize}
\elem

\bpf
Multiplying the eigenequation defining $h$ on the left by $\La^* \Delta_\tau^*$, we obtain
\be
\label{EigenEq}
c \La^* \Delta_\tau^* P_{\cV^\perp} h + d \La^* \Delta_\tau^* \La^* \La h = \la \La^* \Delta_\tau^* h.
\ee
According to \eqref{PropertiesDeltaUnderTF},
we have $\Delta_0^* P_{\cV^\perp} = 0$,
$\Delta_1^* P_{\cV^\perp} = \Delta_1^* - \Delta_0^*$,
$\Delta_1^* \La^* \La = \La$,
and $\Delta_0^* \La^* \La = \Delta_0^* - \Delta_1^* + \La$.
Thus, the relation \eqref{EigenEq} specified to $\tau = 0$ and to $\tau=1$ yields
\begin{align}
\label{EqDelta0h}
d \La^* \Delta_0^* h - d \La^* \Delta_1^* h + d \La^* \La h & = \la \La^* \Delta_0^* h,\\
\label{EqDelta1h}
c \La^* \Delta_1^* h - c \La^* \Delta_0^* h + d \La^* \La h & = \la \La^* \Delta_1^* h. 
\end{align}
Subtracting \eqref{EqDelta1h} from \eqref{EqDelta0h} yields
$(c+d) (\La^* \Delta_0^* h - \La^* \Delta_1^* h) = \la  (\La^* \Delta_0^* h - \La^* \Delta_1^* h)$.
Therefore, we derive that $\La^* \Delta_0^* h = \La^* \Delta_1^* h$ provided $\la \not= c+d$.
In this case,
the equations \eqref{EqDelta0h}-\eqref{EqDelta1h}  reduce to
$\La^* \Delta_0^* h = \La^* \Delta_1^* h = (d / \la) \La^* \La h$.
In view of  $\Delta_\tau = (1-\tau)\Delta_0 + \tau \Delta_1$,
we arrive at $\La^* \Delta_\tau^* h = (d / \la) \La^* \La h$ for any $\tau \in [0,1]$.
The relation $ (\Id  - \La^* \Delta_\tau^*)h= (c/\la) P_{\cV^\perp} h$
follows from the eigenequation rewritten as $(c/\la) P_{\cV^\perp} h + (d/\la) \La^* \La h = h$.

It remains to deal with the case $\la = c+d$.
Notice that this case is not vacuous,
as it is equivalent to $h \in \cV^\perp \cap {\rm ran}(\La^*\La)$,
which is nontrivial by a dimension argument involving assumption \eqref{Assum_VK}.
To see this equivalence, notice that $h \in \cV^\perp \cap {\rm ran}(\La^*\La)$
 clearly implies $c P_{\cV^\perp} h + d  \La^* \La h = (c+d) h$,
while the latter eigenequation forces $ c  \| P_{\cV^\perp}h\|^2 + d  \|\La^* \La h\|^2  = (c+d) \|h\|^2$,  hence $\| P_{\cV^\perp}h\|^2 = \|h\|^2$ and $\|\La^* \La h\|^2 = \|h\|^2$,
i.e., $h \in \cV^\perp$ and $h \in {\rm ran}(\La^*\La)$.
We now consider such an eigenvector~$h$ associated with the eigenvalue $c+d$:
in view of $h \in \cV^\perp \cap {\rm ran}(\La^*\La)$, 
we remark that $\Delta_0^* h = \Delta_0^* P_{\cV^\perp} h =0$
and that $\Delta_1^* h = \Delta_1^* \La^* \La  h = \La h$.
We deduce that $\La^* \Delta_\tau^* h = (1-\tau) \La^* \Delta_0^* h + \tau \La^* \Delta_1^* h = \tau \La^* \La h = \tau h$
and in turn that $(\Id - \La^* \Delta_\tau^*) h = (1-\tau) h$.
\epf

We are now ready to establish the main result of this subsection.

\bpf[Proof of Theorem \ref{ThmIndepTau}]
Let $\tau \in [0,1]$ be fixed throughout.
As announced earlier,
our objective is to establish that,
thanks to $\La \La^* = \Id_{\bR^m}$, the condition $c P_{\cV^\perp} + d \La^* \La \succeq \Id$ implies the condition
$$
\bbmx
c P_{\cV^\perp} & | & 0\\
\hline
0 & | & d \, \Id_{\bR^m}
\ebmx
\succeq \bbmx \Id - \La^* \Delta_\tau^*\\ \hline  \Delta_\tau^* \ebmx
\bbmx \Id - \Delta_\tau \La \; | \; \Delta_\tau \ebmx,
$$
or equivalently the condition
$$
\bbmx
c P_{\cV^\perp} & | & 0\\
\hline
0 & | & d \La^* \La
\ebmx
\succeq \bbmx \Id - \La^* \Delta_\tau^*\\ \hline  \La ^* \Delta_\tau^* \ebmx
\bbmx \Id - \Delta_\tau \La \; | \; \Delta_\tau \La \ebmx.
$$
The equivalence of these conditions is seen as follows:
the former implies the latter by multiplying on the left by 
$\bbmx  \Id & | & 0\\ \hline 0 & | & \La^* \ebmx$
and on the right by $\bbmx  \Id & | & 0\\ \hline 0 & | & \La \ebmx$,
while the latter implies the former
under the assumption $\La \La^* = \Id_{\bR^m}$
by multiplying on the left by 
$\bbmx  \Id & | & 0\\ \hline 0 & | & \La \ebmx$
and on the right by $\bbmx  \Id & | & 0\\ \hline 0 & | & \La^* \ebmx$.
As a matter of fact,
according to a classical result about Schur complements,
see e.g.  \cite[Section A.5.5]{boyd2004convex},
the latter is further equivalent to 
$$
\bbmx
\Id & | &  \Id - \Delta_\tau \La & | & \Delta_\tau \La \\
\hline
\Id - \La^* \Delta_\tau^* & | & c P_{\cV^\perp} & | & 0\\
\hline
\La^* \Delta_\tau^* & | & 0 & | & d \La^* \La
\ebmx
\succeq 0.
$$
Thus, considering $f,g,h \in H$,
our objective is to prove the nonnegativity of the inner product
\begin{align*}
{\rm ip} := \, & 
\left\langle
\bbmx
\Id & | &  \Id - \Delta_\tau \La & | & \Delta_\tau \La \\
\hline
\Id - \La^* \Delta_\tau^* & | & c P_{\cV^\perp} & | & 0\\
\hline
\La^* \Delta_\tau^* & | & 0 & | & d \La^* \La
\ebmx
\bbmx
f \\ \hline g \\ \hline h
\ebmx
, \bbmx
f \\ \hline g \\ \hline h
\ebmx
\right\rangle\\
 = \, & 
\langle f,f \rangle + c \langle P_{\cV^\perp} g, g \rangle + d \langle \La^* \La h, h \rangle
+ 2 \langle (\Id - \La^* \Delta_\tau^*)f,g  \rangle + 2 \langle \La^* \Delta_\tau^* f, h \rangle.  
\end{align*}
Let us decompose $f$, $g$, and $h$ as $f=f'+f''$, $g=g'+g''$, and $h=h'+h''$,
where $f'$, $g'$, and~$h'$ belong to the space $H'$ spanned by eigenvectors of $c P_{\cV^\perp} + d \La^* \La$ corresponding to eigenvalues $\la \not= c+d$
and where $f''$, $g''$, and $h''$ belong to the eigenspace $H''$ of $c P_{\cV^\perp} + d \La^* \La$ corresponding to the eigenvalue $\la= c+d$,
i.e., $H'' = \cV^\perp \cap {\rm ran}(\La^* \La)$.
We take notice of the fact that the spaces $H'$ and $H''$ are orthogonal.
With this decomposition,  the above inner product  becomes
$$
{\rm ip}  = {\rm ip}' + {\rm ip} '' + {\rm ip}''', 
$$
where we have set
\begin{align*}
%\label{IP1}
{\rm ip}'  & =
\langle f',f' \rangle + c \langle P_{\cV^\perp} g', g' \rangle + d \langle \La^* \La h,' h' \rangle
+ 2 \langle (\Id - \La^* \Delta_\tau^*)f',g'  \rangle + 2 \langle \La^* \Delta_\tau^* f', h' \rangle,  \\[5pt] 
%\label{IP2}
{\rm ip}'' & =
\langle f'',f'' \rangle + c \langle P_{\cV^\perp} g'', g'' \rangle + d \langle \La^* \La h'', h'' \rangle
+ 2 \langle (\Id - \La^* \Delta_\tau^*)f'',g''  \rangle + 2 \langle \La^* \Delta_\tau^* f'', h'' \rangle , \\[5pt]
%\label{IP3}
{\rm ip}''' & = 2 \langle f', f'' \rangle + 2 \langle P_{\cV^\perp} g', g'' \rangle + 2 \langle \La^* \La h', h'' \rangle
+ 2 \langle (\Id - \La^* \Delta_\tau^*)f',g''  \rangle + 2 \langle \La^* \Delta_\tau^* f', h'' \rangle\\
\nonumber
&  \qquad \qquad \qquad \qquad \qquad \qquad \qquad \qquad \qquad \,
+ 2 \langle (\Id - \La^* \Delta_\tau^*)f'',g'  \rangle + 2 \langle \La^* \Delta_\tau^* f'', h' \rangle .
\end{align*}
We first remark that the terms in ${\rm ip}'''$ are all zero:
first, it is clear that $\langle f', f'' \rangle = 0$;
then, one has $ \langle P_{\cV^\perp} g', g'' \rangle =  \langle  g',  P_{\cV^\perp}g'' \rangle =  \langle  g',  g'' \rangle =0 $
and $\langle \La^* \La h', h'' \rangle = 0$ is obtained similarly;
next,  Lemma \ref{LemEigcd} ensures that 
$ \langle (\Id - \La^* \Delta_\tau^*)f'',g'  \rangle  =  (1-\tau) \langle f'',g'  \rangle = 0$ and $\langle \La^* \Delta_\tau^* f'', h' \rangle = 0$ is obtained similarly;
last,  writing $f' = \sum_i f_i$ where the $f_i \in H'$ are orthogonal eigenvectors of $c P_{\cV^\perp} + d \La^* \La$ corresponding to eigenvalues $\la_i < c+d$, we derive from Lemma \ref{LemEigcd} that
$$
 \langle (\Id - \La^* \Delta_\tau^*)f',g''  \rangle
 = \sum_i \f{c}{\la_i} \langle P_{\cV^\perp} f_i, g'' \rangle
 = \sum_i \f{c }{\la_i} \langle f_i,  P_{\cV^\perp} g'' \rangle
 = \sum_i \f{c}{\la_i} \langle f_i,  g'' \rangle = 0,
$$
and  $\langle \La^* \Delta_\tau^* f', h'' \rangle = 0$ is obtained similarly.
As a result, we have ${\rm ip}'''=0$.

We now turn to the quantity ${\rm ip}'$.
Exploiting Lemma \ref{LemEigcd} again,  we write
\begin{align*}
{\rm ip}' &=  \langle f' , f' \rangle +  c \langle P_{\cV^\perp} g', g' \rangle + d \langle \La^* \La h,' h' \rangle
+ 2 \bigg\langle \sum_i \f{c}{\la_i} P_{\cV^\perp} f_i,g'  \bigg\rangle 
+ 2 \bigg\langle \sum_i \f{d}{\la_i}  \La^* \La f_i, h' \bigg\rangle\\
&= \langle f' , f' \rangle 
+  c \Bigg(
\langle P_{\cV^\perp} g',   P_{\cV^\perp} g' \rangle
+ 2  \bigg\langle  
\sum_i \f{1}{\la_i} P_{\cV^\perp} f_i,  P_{\cV^\perp} g' \bigg\rangle
\Bigg)\\
& \phantom{= \langle f' , f' \rangle }\;
+ d \Bigg(
\langle \La^* \La h',  \La^* \La  h' \rangle
+ 2 \bigg\langle
\sum_i \f{1}{\la_i} \La^* \La f_i ,  \La^* \La h' \bigg\rangle
\Bigg)\\
& = \langle f' , f' \rangle 
+  c \Bigg( \bigg\| P_{\cV^\perp} g' +  \sum_i \f{1}{\la_i}  P_{\cV^\perp} f_i \bigg\|^2 - \bigg\| \sum_i \f{1}{\la_i}  P_{\cV^\perp} f_i \bigg\|^2 \Bigg)\\
& \phantom{= \langle f' , f' \rangle }\;
+ d \Bigg( \bigg\| \La^* \La h' +  \sum_i \f{1}{\la_i} \La^* \La f_i  \bigg\|^2
- \bigg\| \sum_i \f{1}{\la_i} \La^* \La f_i  \bigg\|^2
 \Bigg).
\end{align*} 
At this point, we can bound $ {\rm ip}'$ from below as
\begin{align*}
 {\rm ip}' & \ge \langle f' , f' \rangle  
 - \Bigg( c \bigg\|  P_{\cV^\perp}  \Big( \sum_i \f{1}{\la_i}  f_i \Big) \bigg\|^2 + d \bigg\| \La^* \La  \Big( \sum_i \f{1}{\la_i} f_i \Big)  \bigg\|^2 \Bigg)\\
 & =  \langle f' , f' \rangle  
 - \bigg\langle \Big( c P_{\cV^\perp} + d \La^* \La \Big) \Big(  \sum_i \f{1}{\la_i}  f_i \Big),  \Big( \sum_i \f{1}{\la_i}  f_i \Big) \bigg\rangle\\
& = \sum_i \|f_i\|^2  
 - \bigg\langle  \sum_i  f_i,  \sum_i \f{1}{\la_i}  f_i \bigg\rangle
 =  \sum_i \|f_i\|^2  \bigg( 1 - \f{1}{\la_i} \bigg).
\end{align*}
This shows that $ {\rm ip}' \ge 0$ since the condition $c P_{\cV^\perp} + d \La^* \La \succeq \Id$ ensures that $\la_i \ge 1$ for every $i$.

Finally,  Lemma \ref{LemEigcd}  also helps us to bound the quantity ${\rm ip}''$ from below according to
\begin{align*}
{\rm ip}''& =  \|f''\|^2 + c \| g'' \|^2 + d \| h''\|^2 + 2 (1-\tau) \langle f'', g'' \rangle + 2 \tau \langle f'',h'' \rangle \\
& = (1-\tau) \big( \|f''\|^2  +  2 \langle f'', g'' \rangle \big)
+ \tau \big(  \|f''\|^2  +  2 \langle f'',h'' \rangle  \big) + c \| g'' \|^2 + d \| h''\|^2 \\
& \ge  - (1-\tau)  \| g'' \|^2  -  \tau  \|h'' \|^2 
+  c \| g'' \|^2 + d \| h''\|^2.
\end{align*}
This allows us to obtain ${\rm ip}'' \ge 0$
since the condition $c P_{\cV^\perp} + d \La^* \La \succeq \Id$ ensures that 
 $c \ge 1 $ and $d \ge 1 $.
Altogether, we have shown that ${\rm ip} = {\rm ip}' + {\rm ip}'' + {\rm ip}''' \ge 0$, which concludes the proof.
\epf

\brk
The value of the minimal global worst-case error can, in general, be computed 
by solving the semidefinite program  \eqref{LB4GWCE} characterizing the lower bound ${\rm lb}$.
In the case where $\La \La^* = \Id_{\bR^m}$,
it can also be computed without resorting to semidefinite programming.
Precisely,
if $\tau_\sharp$ denotes the (unique) $\tau$ between $1/2$ and $\eps/(\eps+\eta)$ such that
\be
\la_{\min}((1-\tau) P_{\cV^\perp} + \tau \La^* \La) =  \f{(1-\tau)^2 \eps^2 - \tau^2 \eta^2}{(1-\tau) \eps^2 - \tau \eta^2 }
\ee
and if $\la_\sharp$ denotes $\la_{\min}((1-\tau_\sharp) P_{\cV^\perp} + \tau_\sharp \La^* \La)$,
then we claim that,  for any $\tau \in [0,1]$, 
$$
{\rm gwce}(\Delta_\tau)^2 = \f{1-\tau_\sharp}{\la_{\sharp}} \eps^2 + \f{\tau_\sharp}{\la_{\sharp}} \eta^2.
$$
Indeed, since we now know that the global worst-case error ${\rm gwce}(\Delta_\tau)$ equals its lower bound ${\rm lb}$ independently of $\tau \in [0,1]$
and since $c_\sharp:=  (1-\tau_\sharp)/\la_{\sharp}$
and $d_\sharp := \tau_\sharp/\la_{\sharp}$
are feasible for the semidefinite program  \eqref{LB4GWCE} characterizing ${\rm lb}$,
we obtain 
\be
\label{4Claim1}
{\rm gwce}(\Delta_\tau)^2 \le \f{1-\tau_\sharp}{\la_{\sharp}} \eps^2 + \f{\tau_\sharp}{\la_{\sharp}} \eta^2.
\ee
Moreover,  going back to the proof of Theorem \ref{ThmLocalTF},
we recognize that the choice of $\tau_\sharp$ here
corresponds to the instance $y=0$ there.
This instance comes with $f_\sharp$ being equal to zero
and with $h_\sharp$ being equal to a properly normalized eigenvector of $(1-\tau_\sharp) P_{\cV^\perp} + \tau_\sharp \La^* \La$ corresponding to the eigenvalue~$\la_\sharp$.
The identities \eqref{Identities4Later} now read $\|P_{\cV^\perp} h_\sharp\|^2 = \eps^2$ and $\|\La^* \La h_\sharp\|^2 = \eta^2$, i.e.,  $\| \La h_\sharp\|^2 = \eta^2$.
Setting $f = h_\sharp$ and $e = -\La h_\sharp$,
which satisfy $\|P_{\cV^\perp} f\|= \eps$
and $\|e\| = \eta$,
the very definition of the global worst-case error yields
\begin{align}
\label{4Claim2}
{\rm gwce}(\Delta_\tau)^2 & \ge \|f - \Delta_\tau(\La f + e)\|^2
%= \|h_\sharp - \Delta_\tau(0)\|^2
= \|h_\sharp\|^2\\
\nonumber
& = \f{1}{\la_\sharp} \big\langle
 \big( (1-\tau_\sharp) P_{\cV^\perp} + \tau_\sharp \La^* \La \big) h_\sharp, h_\sharp
\big\rangle\\
\nonumber
& = \f{1-\tau_\sharp}{\la_\sharp} \|P_{\cV^\perp} h_\sharp\|^2 
+ \f{\tau_\sharp}{\la_\sharp} \|\La^* \La h_\sharp\|^2\\
\nonumber
& =  \f{1-\tau_\sharp}{\la_{\sharp}} \eps^2 + \f{\tau_\sharp}{\la_{\sharp}} \eta^2.
\end{align}
Together,  the inequalities \eqref{4Claim1} and \eqref{4Claim2} justify our claim about the value of the global worst-case error.
In passing, it is worth noticing that the above argument reveals that $f = h_\sharp$ and $e = -\La h_\sharp$ are extremal in the defining expression for the global worst-case error of the regularization map $\Delta_\tau$
independently of the parameter $\tau \in [0,1]$.
\erk

\bibliography{refs}

\newpage

\section*{Appendix}

This additional section collects justifications for a few facts that were mentioned  but not explained in the main text.
These facts are:
the uniqueness of a Chebyshev center for the model- and data-consistent set (see page~\pageref{PageCCUnique}),
the efficient computation of the solution to \eqref{Reg} when $\La \La^* = \Id_{\bR^m}$ (see page~\pageref{PageCompReg}),
the form of Newton method when solving equation \eqref{EqForLocal} (see page~\pageref{PageEigEq}),
and the reason why the constraint of \eqref{GWCETau}
always implies the constraint of \eqref{LB4GWCE} (see pages~\pageref{PageImplies1} and~\pageref{PageImplies2}).

\paragraph{Uniqueness of the Chebyshev center.}

Let $\wh{f_1},\wh{f_2}$ be two Chebyshev centers,
i.e.,  minimizers of $\max\{ \|f-g\|: \|P_{\cV^\perp} g \| \le \eps, \|\La g - y\| \le \eta \}$
and let $\mu$ be the value of the minimum.
Consider $\ol{g} \in H$ such that $\|(\wh{f_1}+\wh{f_2})/2 - \ol{g}\|
= \max\{ \|(\wh{f_1}+\wh{f_2})/2 - g\|  :  \|P_{\cV^\perp} g \| \le \eps, \|\La g - y\| \le \eta \}$.
Then  
\begin{align*}
\mu & \le \|(\wh{f_1}+\wh{f_2})/2 - \ol{g}\|
\le \f{1}{2} \|\wh{f_1} - \ol{g}\| + \f{1}{2} \| \wh{f_2} - \ol{g}\|\\
& \le \f{1}{2} \max\{ \|\wh{f_1}-g\|: \|P_{\cV^\perp} g \| \le \eps, \|\La g - y\| \le \eta \}
+ \f{1}{2} \max\{ \|\wh{f_2}-g\|: \|P_{\cV^\perp} g \| \le \eps, \|\La g - y\| \le \eta \}\\
& = \f{1}{2} \mu + \f{1}{2} \mu = \mu.
\end{align*}
Thus, equality must hold all the way through.
This implies that $\wh{f_1} - \ol{g} = \wh{f_2} - \ol{g}$,
i.e., that $\wh{f_1} = \wh{f_2}$,
as expected.

\paragraph{Computation of the regularized solution.}
Let $(v_1,\ldots,v_n)$ be a basis for $\cV$ 
and let $u_1,\ldots,u_m$ denote the Riesz representers of the observation functionals $\la_1,\ldots,\la_m$,
which form an orthonormal basis for ${\rm ran}(\La^*)$ under the assumption that $\La \La^* = \Id_{\bR^m}$.
With $C \in \bR^{m \times n}$ representing the cross-gramian with entries $\langle u_i,v_j \rangle = \la_i(v_j)$,
the solution to the regularization program \eqref{Reg} is given,
even when $H$ is infinite dimensional,
by
$$
f_\tau = \tau \sum_{i=1}^m a_i u_i +  \sum_{j=1}^n b_j v_j,
$$
where the coefficient vectors $a \in \bR^m$ and $b \in \bR^n$ are computed according to
$$
b = \big( C^\top C \big)^{-1} C^\top y
\qquad \mbox{and} \qquad
a = y-Cb.
$$
This is fairly easy to see for $\tau =0$
and it has been established in \citep*[Theorem 2]{foucart2020learning} for $\tau=1$,
so the general result follows from Proposition \ref{PropLinearExpr}.
Alternatively,  it can be obtained by replicating the steps from the proof of the case $\tau = 1$ with minor changes.

\paragraph{Newton method.}
Equation \eqref{EqForLocal} takes the form $F(\tau) = 0$,
where 
$$
F(\tau) = 
\la_{\min}((1-\tau) R + \tau S) - \f{(1-\tau)^2 \eps^2 - \tau^2 \eta^2}{(1-\tau) \eps^2 - \tau \eta^2 + (1-\tau)\tau(1-2\tau) \delta^2}.
$$
Newton method produces a sequence $(\tau_k)_{k \ge 0}$ converging to a solution using the recursion
\be
\label{NewtonM}
\tau_{k+1} = \tau_k - \f{F(\tau_k)}{F'(\tau_k)},
\qquad k \ge 0.
\ee
In order to apply this method,
we need the ability to compute the derivative of $F$ with respect to $\tau$.
Setting $\la_{\min} = \la_{\min}((1-\tau)R+\tau S)$,
this essentially reduces to the computation of $d \la_{\min}/d\tau$,
which is performed via the argument below.
Note that the argument is not rigorous,
as we take for granted the differentiability of the eigenvalue $\la_{\min}$ and of a normalized eigenvector $h$ associated with~it.
However, nothing prevents us from applying the scheme \eqref{NewtonM}
using the expression for $d \la_{\min}/d\tau$ given in \eqref{D4Newton} below
and agree that a solution has been found if the output $\tau_K$ satisfies $F(\tau_K) < \iota$ for some prescribed tolerance $\iota >0$.
Now, the argument
starts form the identities 
$$
((1-\tau)R+\tau S) h = \la_{\min} h
\qquad \mbox{and} \qquad 
\langle h,h \rangle =1,
$$
which we differentiate to obtain
$$
(S-R)h + ((1-\tau)R+\tau S) \f{dh}{d\tau}
 = \f{d \la_{\min}}{d \tau} h + \la_{\min} \f{dh}{d \tau}
 \qquad \mbox{and} \qquad 
 2 \Big\langle h, \f{dh}{d\tau} \Big\rangle = 0.
$$
By taking the inner product with $h$ in the first identity and using the second identity,
we derive
$$
\langle (S-R)h , h \rangle = \f{d \la_{\min}}{d\tau}, 
\qquad \mbox{i.e., }\qquad
\f{d \la_{\min}}{d\tau} = \|S h\|^2 - \|Rh\|^2.
$$
According to Lemma \ref{LemEigenStuff},
this expression can be transformed, after some work, into 
\be
\label{D4Newton}
\f{d \la_{\min}}{d\tau} = \f{1-2\tau}{\tau(1-\tau)} \, \f{\la_{\min}(1-\la_{\min})}{1-2\la_{\min}}.
\ee

\paragraph{Relation between semidefinite constraints.}

Suppose that the constraint of \eqref{GWCETau} holds for a regularization map $\Delta_\tau$.
In view of the expressions 
$$
\Delta_\tau = \big( (1-\tau) P_{\cV^\perp} + \tau \La^* \La \big)^{-1} (\tau \La^*)
\quad \mbox{and} \quad
\Id - \Delta_\tau \La = \big( (1-\tau) P_{\cV^\perp} + \tau \La^* \La \big)^{-1} ((1-\tau) P_{\cV^\perp}),
$$ 
this constraint also reads
$$
\bbmx
c P_{\cV^\perp} & | & 0\\
\hline
0 & | & d \,  \Id_{\bR^m}
\ebmx 
\succeq
 \bbmx (1-\tau) P_{\cV^\perp}\\ \hline  \tau \La \ebmx
 \big( (1-\tau) P_{\cV^\perp} + \tau \La^* \La \big)^{-2}
\bbmx (1-\tau) P_{\cV^\perp} \; | \; \tau \La^* \ebmx. 
$$
Multiplying on the left by $\bbmx P_{\cV^\perp} \; | \; \La^* \ebmx$
and on the right by $\bbmx P_{\cV^\perp} \\ \hline \La \ebmx$ yields
$$
c P_{\cV^\perp}  + d \La^* \La \succeq
 ( (1-\tau) P_{\cV^\perp} + \tau \La^* \La )
 \big( (1-\tau) P_{\cV^\perp} + \tau \La^* \La \big)^{-2}
 ( (1-\tau) P_{\cV^\perp} + \tau \La^* \La )
 = \Id.
$$
This is the constraint of \eqref{LB4GWCE}.

\end{document}